\documentclass{amsart}
\usepackage{graphicx,amscd,color,amsmath,amsfonts,amssymb,geometry}
\providecommand{\U}[1]{\protect\rule{.1in}{.1in}}

\newtheorem{theorem}{Theorem}[section]
\newtheorem{proposition}[theorem]{Proposition}

\newtheorem{lemma}[theorem]{Lemma}

\newtheorem{definition}[theorem]{Definition}
\numberwithin{equation}{section}

\begin{document}
\author[A. T. Bernardino \and D. Pellegrino \and J.B. Seoane-Sep\'{u}lveda \and M.L.V. Souza]{A.T. Bernardino \and D. Pellegrino\textsuperscript{*} \and J.B. Seoane-Sep\'{u}lveda\textsuperscript{**} \and M.L.V. Souza}
\address{Centro de Ensino Superior do Serid\'{o},\newline\indent Universidade Federal do Rio Grande do Norte, \newline\indent Rua Joaquim Greg\'{o}rio, S/N - Penedo, \newline\indent Caic\'{o}, 59300-000, Brazil.}
\email{thiagobernardino@yahoo.com.br}
\address{Departamento de Matem\'{a}tica,\newline\indent Universidade Federal da Para\'{\i}ba,\newline\indent 58.051-900 - Jo\~{a}o Pessoa, Brazil.}
\email{pellegrino@pq.cnpq.br}
\thanks{\textsuperscript{*}Supported by CNPq Grant 301237/2009-3.}
\address{Departamento de An\'{a}lisis Matem\'{a}tico,\newline\indent Facultad de Ciencias Matem\'{a}ticas, \newline\indent Plaza de Ciencias 3, \newline\indent Universidad Complutense de Madrid,\newline\indent Madrid, 28040, Spain.}
\email{jseoane@mat.ucm.es}
\thanks{\textsuperscript{**}Supported by the Spanish Ministry of Science and Innovation, grant MTM2009-07848.}
\address{Departamento de Matem\'{a}tica/ICENE,\newline\indent UFTM - Universidade Federal do Tri\^angulo Mineiro,\newline\indent
Rua Get\'ulio Guarit\'a, 159, \newline\indent CEP 38.025-440 - Uberaba-MG, Brazil.}
\email{marcelalvsouza@gmail.com}
\title[Absolutely summing operators revisited]{Absolutely summing operators revisited: new directions in the nonlinear theory}
\begin{abstract}
In the last decades many authors have become interested in the study of multilinear and polynomial generalizations of families of operator ideals (such as, for instance, the ideal of absolutely summing operators). However, these generalizations must keep the essence of the given operator ideal and there seems not to be a universal method to achieve this. The main task of this paper is to discuss, study, and introduce multilinear and polynomial extensions of the aforementioned operator ideals taking into account the already existing methods of evaluating the adequacy of such generalizations. Besides this subject's intrinsic mathematical interest, the main motivation is our belief (based on facts that shall be presented) that some of the already existing approaches are not adequate.
\end{abstract}
\maketitle


\section{Introduction and historical background}

A well-known fact from an undergraduate Analysis course states that, in
$\mathbb{R}$, a series converges absolutely if and only if it is
unconditionally convergent; this result was proved by J.P.G.L. Dirichlet in
1829. For infinite-dimensional Banach spaces the situation is quite different:
on the one hand for $\ell_{p}$ spaces with $1<p<\infty$, for example, it is
quite easy to construct an unconditionally convergent series which fails to be
absolutely convergent. On the other hand, for $\ell_{1}$ and some other Banach
spaces the answer to this problem is far from being straightforward. The
special case of $\ell_{1}$ was solved in 1947 by M.S. Macphail \cite{Mac}
through a very elaborated construction.

The question of whether every infinite-dimensional Banach space has an
unconditionally convergent series which fails to be absolutely convergent was
raised by Banach \cite[p. 40]{Banach32} (see also Problem 122 in the Scottish
Book \cite{Mau}, proposed by S. Mazur and W. Orlicz). In 1950, A. Dvoretzky
and C.A. Rogers \cite{DR} solved this question in the positive:

\bigskip

\textbf{Theorem} (Dvoretzky-Rogers, 1950). The unconditionally convergent
series and absolutely summing convergent series coincide in a Banach space $E
$ if and only if $\dim E=\infty.$

\bigskip

The above result encouraged the curiosity of the genius of A. Grothendieck,
who rapidly presented a different proof of this result in his Ph.D.
dissertation \cite{Gro1955}. Grothendieck's famous R\'{e}sum\'{e} \cite{gro}
(see also \cite{dies01} for a modern and thorough study) and \cite{Gro1955}
are, essentially, the beginning of the theory of absolutely summing operators.
More precisely, in view of Dvoretzky-Rogers' striking result, the idea of
investigating linear operators that transform unconditionally convergent
series into absolutely convergent series seemed natural and was the birth of
the notion of absolutely summing operators (a linear operator $u:E\rightarrow
F$ is absolutely summing if ${\textstyle\sum} u(x_{j})$ is absolutely
convergent whenever ${\textstyle\sum} x_{j}$ is unconditionally convergent).
Soon after, Grothendieck proved a quite surprising result asserting that every
continuous linear operator from $\ell_{1}$ to $\ell_{2}$ (or to any Hilbert
space) is absolutely summing (this kind of result is now called a
\emph{coincidence theorem}). This result is a consequence of an intriguing
inequality which Grothendieck himself called \textquotedblleft the fundamental
theorem of the metric theory of tensor products\textquotedblright.
Grothendieck's inequality has important applications (\cite{AAA, FFF}) and
still has some hidden mysteries such as the precise value of Grothendieck's
constant. For a recent work on the estimates for Grothendieck's constant we
refer to \cite{NN}.

The modern notion of absolutely $(p;q)$-summing operators was introduced in
the 1960's by A. Pietsch \cite{stu} and B. Mitiagin and A. Pe\l czy\'{n}ski
\cite{MPel}. Besides its intrinsic mathematical interest and deep mathematical
motivation, it has shown to be a very important tool in general Banach space
theory. For instance, and just to cite some, using the theory of absolutely
summing operators one can show that every normalized unconditional basis of
$\ell_{1}$ is equivalent to the unit vector basis of $\ell_{1}$ and also that,
for $1<p<\infty$, there is a normalized unconditional basis of $\ell_{p}$
which is not equivalent to the unit vector basis of $\ell_{p}$.

Throughout this paper $\mathbb{N}$ represents the set of all positive integers
and $\mathbb{N}_{m}:=\{1,...,m\}$. Also, $E,E_{1},\ldots,E_{n},F,G,G_{1}%
,...,G_{n},H$ will stand for Banach spaces over $\mathbb{K}=\mathbb{R}$ or
$\mathbb{C}$, the topological dual of $E$ is represented by $E^{\ast}$ and
$B_{E^{\ast}}$ denotes its closed unit ball. The symbol $W\left(  B_{E^{\ast}%
}\right)  $ represents the probability measures in the Borel sets of
$B_{E^{\ast}}$ with the weak-star topology. We will denote the space of all
continuous $n$-linear operators from $E_{1}\times\cdots\times E_{n}$ into $F$
by $\mathcal{L}(E_{1},\ldots,E_{n};F) $ or $\mathcal{L}_{n}(E_{1},\ldots
,E_{n};F).$ Also, we recall that an $n$-homogeneous polynomial $P:E\rightarrow
F$ is a map so that $P(x)=\check{P}(x,\ldots,x),$ where $\check{P}$ represents
the unique symmetric $n$-linear map associated to $P$. The corresponding space
(endowed with the sup norm) is represented by $\mathcal{P}(^{n}E;F)$. For the
theory of polynomials and multilinear operators acting on Banach spaces we
refer to \cite{Dineen, Mu}.

For $0<p<\infty$, the space of all sequences $\left(  x_{j}\right)
_{j=1}^{\infty}$ in $E$ such that $\left(  \varphi\left(  x_{j}\right)
\right)  _{j=1}^{\infty}\in\ell_{p}$, for every $\varphi\in E^{\ast}$ is
denoted by $\ell_{p}^{w}\left(  E\right)  .$ When endowed with the norm
($p$-norm if $0<p<1 $)
\[
\left\Vert \left(  x_{j}\right)  _{j=1}^{\infty}\right\Vert _{w,p}%
{\small :=}\sup\{({\textstyle\sum\limits_{j=1}^{\infty}}\left\vert
\varphi\left(  x_{j}\right)  \right\vert ^{p})^{1/p}:\varphi\in B_{E^{\ast}%
}\}{\small ,}%
\]
the space $\ell_{p}^{w}\left(  E\right)  $ is complete. We recall that if
$0<q\leq p<\infty$ a continuous linear operator $u:E\rightarrow F$ is
absolutely $(p;q)$-summing if $\left(  u(x_{j})\right)  _{j=1}^{\infty}\in
\ell_{p}\left(  F\right)  $ whenever $\left(  x_{j}\right)  _{j=1}^{\infty}%
\in\ell_{q}^{w}\left(  E\right)  .$ In this case we write $u\in\Pi
_{(p;q)}(E;F)$. For $p=q=1$ this notion coincides with the concept of
absolutely summing operator. For classical results on absolutely summing
operators we refer to \cite{dies0, pisier, T1} and references therein (recent
results can also be checked in \cite{PellZ, Lima, Ok}). The concept of
absolutely summing operators has some natural linear extensions such as the
notions of mixing $(p;q)$-summing operators (due to A. Pietsch and B. Maurey)
and $(p;q;r)$-summing operators (due to A. Pietsch). It is worth mentioning
that these concepts were not just constructed to simply generalize the notion
of absolutely $(p;q)$-summing operators; these notions have their particular
reasons to be investigated (see \cite[p. 359]{hist}).

In the 1980's, Pietsch \cite{Pie} suggested a multilinear approach to the
theory of absolutely summing operators and, more generally, to the theory of
operator ideals. Since then, several authors were attracted by the subject and
also non-multilinear approaches have appeared (see \cite{CDo, Chen, Junek,
Nach00, MST, adv}). The adequate way of lifting the notion of a given operator
ideal to the multilinear and polynomial settings is a delicate matter. For
example, in the case of the ideal of absolutely summing linear operators,
there are several different approaches to the polynomial and multilinear
contexts (see \cite{adv, ppaa} and references therein). The abstract notions
of (global) holomorphy types (see \cite{BBJMs, Nachbin}), coherent and
compatible ideals (see \cite{CDM09}) shed some light on what kind of approach
is more adequate.

Recently, in 2003, the notion of multiple summing multilinear operators (and
polynomials) was introduced (see \cite{Matos, pv}) but, as a matter of fact,
the origin of this notion dates back to \cite{bh, lit, Ram}. Several
indicators from the theory of summing operators and from the theory of
(multi-) ideals show that this is one of the most adequate approaches to the
nonlinear theory of absolutely summing operators. For results on multiple
summing multilinear operators we refer to \cite{Na, df, davidstudia, pv,
PopJM, Popa}.

Notwithstanding the quick success of the theory of multiple summing
multilinear operators, some recent papers related to multilinear summability
seem to have overlooked its advantages. More precisely, the multilinear
notions of mixing summing operators and absolutely $\left(  p;q;r\right)
$-summing multilinear operators were introduced following a different
perspective (see \cite{achour, Carlos Alberto}). The point is that these
approaches do not carry out the essence of the respective linear concepts and
this lack is clearly corroborated by the notions of coherence, compatibility
and holomorphy types.

In this paper we present multilinear and polynomial notions of absolutely
$\left(  p;q;r\right)  $-summing operators and mixing summing operators which
follow the philosophy of the idea of multiple summability. Among other
results, the adequacy of our approach is evaluated by proving that our new
definitions provide coherent sequences, compatible and also (global)
holomorphy types.

Below we recall the notions of mixing summing operators and absolutely
$\left(  p;q;r\right)  $-summing operators.

\subsection{Mixing summing operators\label{xxzzz}}

Let $0<p\leq s\leq\infty$ and $r$ such that $\frac{1}{r}+\frac{1}{s}=\frac
{1}{p}.$ A sequence $(x_{i})_{i=1}^{\infty}$ in $E$ is $(s;p)$-mixed summable
if
\[
x_{i}=\tau_{i}y_{i}%
\]
with $(\tau_{i})_{i=1}^{\infty}\in\ell_{r}$ and $(y_{i})_{i=1}^{\infty}\in
\ell_{s}^{w}(E)$.

In this case, consider
\[
\left\Vert \left(  x_{i}\right)  _{i=1}^{\infty}\right\Vert _{mx(s,p)}%
:=\inf\left\{  \left\Vert \left(  \tau_{i}\right)  _{i=1}^{\infty}\right\Vert
_{r}\left\Vert \left(  y_{i}\right)  _{i=1}^{\infty}\right\Vert _{w,s}%
\right\}  ,
\]
where the infimum is taken over all possible representations of $\left(
x_{i}\right)  _{i=1}^{\infty}$ in the above form. The space of all
$(s;p)$-mixed summable sequences in $E$ is represented by $\ell_{(s,p)}%
^{mx}(E).$ It is not difficult to prove that $\ell_{(s,p)}^{mx}(E)$ is a
complete normed ($p$-normed if $0<p<1$) space.

It is immediate that, for $0<p\leq s\leq\infty,$ one always has

\begin{itemize}
\item $\ell_{p}(E)\subset\ell_{(s,p)}^{mx}(E)\subset\ell_{p}^{w}(E)$ with%
\begin{equation}
\left\Vert \left(  z_{j}\right)  _{j=1}^{\infty}\right\Vert _{w,p}%
\leq\left\Vert \left(  z_{j}\right)  _{j=1}^{\infty}\right\Vert _{mx(s,p)}%
\leq\left\Vert \left(  z_{j}\right)  _{j=1}^{\infty}\right\Vert _{p},
\label{prop}%
\end{equation}

\item $\ell_{p}^{w}(E)=\ell_{(p,p)}^{mx}(E)$ and $\ell_{p}(E)=\ell
_{(\infty,p)}^{mx}(E)$ isometrically.
\end{itemize}

Let us now recall the linear concept of mixing summing linear operators (see
\cite{pp1}):

Let $0<p\leq s\leq\infty.$ A continuous linear operator $u:E\rightarrow F$ is
mixing $(s,p)$-summing ($u\in\Pi_{mx(s,p)}(E;F)$) if there exists a constant
$\sigma\geq0$ such that%
\begin{equation}
\left\Vert \left(  u(x_{j})\right)  _{j=1}^{m}\right\Vert _{mx(s,p)}\leq
\sigma\left\Vert (x_{j})_{j=1}^{m}\right\Vert _{w,p} \label{rst}%
\end{equation}
for all $x_{1},\ldots,x_{m}\in E$ and $m\in\mathbb{N}.$ The infimum of all
such constants $\sigma$ is represented by $\pi_{mx(s,p)}(u).$

The terminology \textquotedblleft mixing\textquotedblright\ is motivated by
the fact that a continuous linear operator $u:E\rightarrow F$ is $\left(
s,p\right)  $-mixing summing precisely when $u$ maps every weakly $p$-summable
sequence $\left(  x_{i}\right)  _{i=1}^{\infty}$ in $E$ into a sequence which
can be written as a product $\left(  \tau_{i}y_{i}\right)  _{i=1}^{\infty}$ of
an absolutely $r$-summable scalar sequence $\left(  \tau_{i}\right)
_{i=1}^{\infty}$ and a weakly $s$-summable sequence $\left(  y_{i}\right)
_{i=1}^{\infty}$ in $F$, where $\frac{1}{s}+\frac{1}{r}=\frac{1}{p}.$ Many of
the classical results of mixing summing operators are due to B. Maurey
\cite{Maurey} and the theory has shown to be sufficiently rich to be
investigated by its own (see \cite[Section 32]{Flore}).

\subsection{Absolutely $(p;q;r)$-summing operators\label{Sub2}}

The concept of absolutely $(p;q;r)$-summing linear operators is due to A.
Pietsch \cite{pp0, pp1}. If $0<p,q<\infty$ and $0<r\leq\infty$ and%
\[
\frac{1}{p}\leq\frac{1}{q}+\frac{1}{r},
\]
a continuous linear operator $u:E\rightarrow F$ is absolutely $\left(
p;q;r\right)  $-summing ($u\in\Pi_{as\left(  p;q;r\right)  }\left(
E;F\right)  $) if there is a constant $C>0$ such that%
\begin{equation}
\left(  \sum_{j=1}^{m}\left\vert \varphi_{j}\left(  u\left(  x_{j}\right)
\right)  \right\vert ^{p}\right)  ^{\frac{1}{p}}\leq C\left\Vert \left(
x_{j}\right)  _{j=1}^{m}\right\Vert _{w,q}\left\Vert \left(  \varphi
_{j}\right)  _{j=1}^{m}\right\Vert _{w,r} \label{1807}%
\end{equation}
for all positive integer $m$, and all $x_{1},\ldots,x_{m}$ in $E$ and
$\varphi_{1},\ldots,\varphi_{m}$ in $F^{\ast}$. When $r=\infty$, we recover
the classical notion of absolutely $(p;q)$-summing operators.\ For details we
refer to \cite{Lap, pp1, hist}.

The space composed by all continuous linear operators from $E$ to $F$ that are
absolutely $\left(  p;q;r\right)  $-summing shall be represented by
$\Pi_{as\left(  p;q;r\right)  }\left(  E;F\right)  $. The infimum of the
constants $C$ satisfying the inequality (\ref{1807}) defines a norm ($p$-norm
if $0<p<1$) in $\Pi_{as\left(  p;q;r\right)  }\left(  E;F\right)  ,$ denoted
by $\pi_{\left(  p;q;r\right)  }(u).$ If $r=\infty$ we use the classical
notation of absolutely $\left(  p;q\right)  $-summing operators, $\Pi_{\left(
p;q\right)  }\left(  E;F\right)  $ and $\pi_{\left(  p;q\right)  }$ for the norm.

If we allow $\frac{1}{p}>\frac{1}{q}+\frac{1}{r}$ we would have $\Pi
_{as\left(  p;q;r\right)  }\left(  E;F\right)  =\left\{  0\right\}  $ (see
\cite[p. 196]{djt}) and, for this reason, we ask for $\frac{1}{p}\leq\frac
{1}{q}+\frac{1}{r}$ in the definition above.

\subsection{Operator ideals, multi-ideals and polynomial ideals}

The theory of operator ideals goes back to J.W. Calkin \cite{cal}, H. Weyl
\cite{we} and further work of A.\ Grothendieck \cite{grote}. However, only in
the 70's, with A. Pietsch \cite{pp1}, the theory was organized in the modern
presentation (see also \cite{ddjj, HP}). For historical details we suggest
\cite{hist} and for applications we refer to \cite{ddjj}.

An operator ideal $\mathcal{I}$ is a subclass of the class $\mathcal{L}_{1}$
of all continuous linear operators between Banach spaces such that for all
Banach spaces $E$ and $F$ its components%
\[
\mathcal{I}(E;F):=\mathcal{L}_{1}(E;F)\cap\mathcal{I}%
\]
satisfy the following:

(Oa) $\mathcal{I}(E;F)$ is a linear subspace of $\mathcal{L}_{1}(E;F)$ which
contains the finite rank operators.

(Ob) If $u\in\mathcal{I}(E;F)$, $v\in\mathcal{L}_{1}(G;E)$ and $w\in
\mathcal{L}_{1}(F;H)$, then $w\circ u\circ v\in\mathcal{I}(G;H)$.

The operator ideal is called a normed operator ideal if there is a function
$\Vert\cdot\Vert_{\mathcal{I}}\colon\mathcal{I}\longrightarrow\lbrack
0,\infty)$ satisfying

\bigskip

(Ob1) $\Vert\cdot\Vert_{\mathcal{I}}$ restricted to $\mathcal{I}(E;F)$ is a
norm, for all Banach spaces $E$, $F$.

(Ob2) $\Vert P_{1}\colon\mathbb{K}\longrightarrow\mathbb{K}:P_{1}%
(\lambda)=\lambda\Vert_{\mathcal{I}}=1.$

(Ob3) If $u\in\mathcal{I}(E;F)$, $v\in\mathcal{L}_{1}(G;E)$ and $w\in
\mathcal{L}_{1}(F;H)$, then
\[
\Vert w\circ u\circ v\Vert_{\mathcal{I}}\leq\Vert w\Vert\Vert u\Vert
_{\mathcal{I}}\Vert v\Vert.
\]
When $\mathcal{I}(E;F)$ with the norm above is always complete, $\mathcal{I}$
is called a Banach operator ideal.

Absolutely summing operators and the two related aforementioned concepts are
examples of operator ideals. Other examples include the compact, weakly
compact, strictly singular operators, etc.

The notion of multi-ideals is also due to Pietsch \cite{Pie}. For each
positive integer $n$, let $\mathcal{L}_{n}$ denote the class of all continuous
$n$-linear operators between Banach spaces. An ideal of multilinear mappings
(or multi-ideal) $\mathcal{M}$ is a subclass of the class $\mathcal{L}=%
{\textstyle\bigcup\limits_{n=1}^{\infty}}
\mathcal{L}_{n}$ of all continuous multilinear operators between Banach spaces
such that for a positive integer $n$, Banach spaces $E_{1},\ldots,E_{n}$ and
$F$, the components
\[
\mathcal{M}_{n}(E_{1},\ldots,E_{n};F):=\mathcal{L}_{n}(E_{1},\ldots
,E_{n};F)\cap\mathcal{M}%
\]
satisfy:

\bigskip

(Ma) $\mathcal{M}_{n}(E_{1},\ldots,E_{n};F)$ is a linear subspace of
$\mathcal{L}_{n}(E_{1},\ldots,E_{n};F)$ which contains the $n$-linear mappings
of finite type.

(Mb) If $T\in\mathcal{M}_{n}(E_{1},\ldots,E_{n};F)$, $u_{j}\in\mathcal{L}%
_{1}(G_{j};E_{j})$ for $j=1,\ldots,n$ and $v\in\mathcal{L}_{1}(F;H)$, then
\[
v\circ T\circ(u_{1},\ldots,u_{n})\in\mathcal{M}_{n}(G_{1},\ldots,G_{n};H).
\]

Moreover, $\mathcal{M}$ is a (quasi-) normed multi-ideal if there is a
function $\Vert\cdot\Vert_{\mathcal{M}}\colon\mathcal{M}\longrightarrow
\lbrack0,\infty)$ satisfying

\bigskip

(Mb1) $\Vert\cdot\Vert_{\mathcal{M}}$ restricted to $\mathcal{M}_{n}%
(E_{1},\ldots,E_{n};F)$ is a (quasi-) norm, for all Banach spaces
$E_{1},\ldots,E_{n}$ and $F.$

(Mb2) $\Vert T_{n}\colon\mathbb{K}^{n}\longrightarrow\mathbb{K}:T_{n}%
(\lambda_{1},\ldots,\lambda_{n})=\lambda_{1}\cdots\lambda_{n}\Vert
_{\mathcal{M}}=1$ for all $n$,

(Mb3) If $T\in\mathcal{M}_{n}(E_{1},\ldots,E_{n};F)$, $u_{j}\in\mathcal{L}%
_{1}(G_{j};E_{j})$ for $j=1,\ldots,n$ and $v\in\mathcal{L}_{1}(F;H)$, then
\[
\Vert v\circ T\circ(u_{1},\ldots,u_{n})\Vert_{\mathcal{M}}\leq\Vert
v\Vert\Vert T\Vert_{\mathcal{M}}\Vert u_{1}\Vert\cdots\Vert u_{n}\Vert.
\]
When all the components $\mathcal{M}_{n}(E_{1},\ldots,E_{n};F)$ are complete
under this (quasi-) norm, $\mathcal{M}$ is called a (quasi-) Banach
multi-ideal. For a fixed multi-ideal $\mathcal{M}$ and a positive integer $n$,
the class%
\[
\mathcal{M}_{n}:=\cup_{E_{1},\ldots,E_{n},F}\mathcal{M}_{n}\left(  E_{1}
,\ldots,E_{n};F\right)
\]
is called ideal of $n$-linear mappings.

Similarly, for each positive integer $n$, let $\mathcal{P}_{n}$ denote the
class of all continuous $n$-homogeneous polynomials between Banach spaces. A
polynomial ideal $\mathcal{Q}$ is a subclass of the class $\mathcal{P}=
{\textstyle\bigcup\limits_{n=1}^{\infty}} \mathcal{P}_{n}$ of all continuous
homogeneous polynomials between Banach spaces so that for all $n\in\mathbb{N}
$ and all Banach spaces $E$ and $F$, the components
\[
\mathcal{Q}_{n}\left(  ^{n}E;F\right)  :=\mathcal{P}_{n}\left(  ^{n}
E;F\right)  \cap\mathcal{Q}%
\]
satisfy:

(Pa) $\mathcal{Q}_{n}\left(  ^{n}E;F\right)  $ is a linear subspace of
$\mathcal{P}_{n}\left(  ^{n}E;F\right)  $ which contains the finite-type polynomials.

(Pb) If $u\in\mathcal{L}_{1}\left(  G;E\right)  $, $P\in\mathcal{Q}_{n}\left(
^{n}E;F\right)  $ and $w\in\mathcal{L}_{1}\left(  F;H\right)  $, then
\[
w\circ P\circ u\in\mathcal{Q}_{n}\left(  ^{n}G;H\right)  .
\]
\bigskip

If there exists a map $\left\Vert \cdot\right\Vert _{\mathcal{Q}}
:\mathcal{Q}\rightarrow\lbrack0,\infty\lbrack$ satisfying

(Pb1) $\left\Vert \cdot\right\Vert _{\mathcal{Q}}$ restricted to
$\mathcal{Q}_{n}(^{n}E;F)$ is a (quasi-) norm for all Banach spaces $E$ and
$F$ and all $n$;

(Pb2) $\left\Vert P_{n}:\mathbb{K}\rightarrow\mathbb{K};\text{ }P_{n}\left(
\lambda\right)  =\lambda^{n}\right\Vert _{\mathcal{Q}}=1$ for all $n$;

(Pb3) If $u\in\mathcal{L}_{1}(G;E)$, $P\in\mathcal{Q}_{n}(^{n}E;F)$ and
$w\in\mathcal{L}_{1}(F;H),$ then
\[
\left\Vert w\circ P\circ u\right\Vert _{\mathcal{Q}}\leq\left\Vert
w\right\Vert \left\Vert P\right\Vert _{\mathcal{Q}}\left\Vert u\right\Vert
^{n},
\]
$\mathcal{Q}$ is called (quasi-) normed polynomial ideal. If all components
$\mathcal{Q}_{n}\left(  ^{n}E;F\right)  $ are complete, $\left(
\mathcal{Q},\left\Vert \cdot\right\Vert _{\mathcal{Q}}\right)  $ is called a
(quasi-) Banach ideal of polynomials (or (quasi-) Banach polynomial ideal).
For a fixed ideal of polynomials $\mathcal{Q}$ and $n\in\mathbb{N}$, the
class
\[
\mathcal{Q}_{n}:=\cup_{E,F}\mathcal{Q}_{n}\left(  ^{n}E;F\right)
\]
is called ideal of $n$-homogeneous polynomials.

A crucial question in the theory of Banach polynomial ideals (and
multi-ideals) is the following:

\begin{quote}
\textit{Given an operator ideal, is there a natural method to define a related
multi-ideal and polynomial ideal without loosing its essence?}
\end{quote}

As mentioned before, in general a given operator ideal has several different
possible extensions to multi-ideals and polynomial ideals. In an attempt of
filtering what approaches are better than others the notions of coherence,
compatibility (and in some sense holomorphy types) are quite helpful.

In the last decades several authors have been interested in investigating
multilinear and polynomial generalizations of certain operator ideals, such as
the ideal of absolutely summing operators. But the search for the correct
approach is not an easy task. The generalizations must keep the essence of the
given operator ideal and there seems to be no universal receipt for it.

The main goal of this paper is to discuss and introduce multilinear and
polynomial extensions of the aforementioned operator ideals (from Subsections
\ref{xxzzz} and \ref{Sub2}) taking into account the existent methods of
evaluating the adequacy of such generalizations. Besides the intrinsic
mathematical interest of the subject, the main motivation of this paper is
that we believe (based on concrete facts) that the previous approaches were
not adequate.

\section{Coherence and compatibility}

The notions of coherent sequences of ideals of polynomials and compatible
ideals of polynomials, which we recall below, are important tools for
evaluating polynomial extensions of a given operator ideal. The essence of
these concepts rests in the searching of harmony between the levels of
homogeneity ($n$-linearity) of a polynomial ideal and connections
(compatibility) with the case of linear operators ($n=1$). In the following if
$P\in\mathcal{P}\left(  ^{n}E;F\right)  $, then $P_{a^{k}}\in\mathcal{P}%
\left(  ^{n-k}E;F\right)  $ is defined by%
\[
P_{a^{k}}(x):=\check{P}(a,\ldots,a,x,\ldots,x).
\]

\begin{definition}
[Compatible ideals, \cite{CDM09}]\label{IdeaisCompativeis}Let $\mathcal{U}$ be
a normed ideal of linear operators. A normed ideal of $n$-homogeneous
polynomials $\mathcal{U}_{n}$ is compatible with $\mathcal{U}$ if there exist
positive constants $\alpha_{1}$ and $\alpha_{2}$ such that for every Banach
spaces $E$ and $F$, the following conditions hold:

$\left(  i\right)  $ For each $P\in\mathcal{U}_{n}\left(  E;F\right)  $ and
$a\in E$, $P_{a^{n-1}}$ belongs to $\mathcal{U}\left(  E;F\right)  $ and
\[
\left\Vert P_{a^{n-1}}\right\Vert _{\mathcal{U}\left(  E;F\right)  }\leq
\alpha_{1}\left\Vert P\right\Vert _{\mathcal{U}_{n}\left(  E;F\right)
}\left\Vert a\right\Vert ^{n-1}.
\]

$\left(  ii\right)  $ For each $T\in\mathcal{U}\left(  E;F\right)  $ and
$\gamma\in E^{\ast}$, $\gamma^{n-1}T$ belongs to $\mathcal{U}_{n}\left(
E;F\right)  $ and
\[
\left\Vert \gamma^{n-1}T\right\Vert _{\mathcal{U}_{n}\left(  E;F\right)  }%
\leq\alpha_{2}\left\Vert \gamma\right\Vert ^{n-1}\left\Vert T\right\Vert
_{\mathcal{U}\left(  E;F\right)  }.
\]

\end{definition}

For the sake of simplicity, we will sometimes write \textquotedblleft the
sequence $\left(  \mathcal{U}_{n}\right)  _{n=1}^{\infty}$ is compatible with
$\mathcal{U}$\textquotedblright\ instead of writing \textquotedblleft%
$\mathcal{U}_{n}$ is compatible with $\mathcal{U}$ for every $n$%
\textquotedblright. Besides, when we write \textquotedblleft the sequence
$\left(  \mathcal{U}_{n}\right)  _{n=1}^{\infty}$ fails to be compatible with
$\mathcal{U}$\textquotedblright\ we are saying that at least for some $n$, the
ideal $\mathcal{U}_{n}$ is not compatible with $\mathcal{U}$.

\begin{definition}
[Coherent sequence of polynomial ideals \cite{CDM09}]\label{IdeaisCoerentes}%
Consider the sequence $\left(  \mathcal{U}_{k}\right)  _{k=1}^{N}$, where for
each $k$, $\mathcal{U}_{k}$ is an ideal of $k$-homogeneous polynomials and $N$
is eventually infinite. The sequence $\left(  \mathcal{U}_{k}\right)
_{k=1}^{N}$ is a coherent sequence of polynomial ideals if there exist
positive constants $\beta_{1}$ and $\beta_{2}$ such that for every Banach
spaces $E$ and $F$, the following conditions hold for $k\in\{1,\ldots,N-1\}$:

$\left(  i\right)  $ For each $P\in\mathcal{U}_{k+1}\left(  E;F\right)  $ and
$a\in E$, $P_{a}$ belongs to $\mathcal{U}_{k}\left(  E;F\right)  $ and
\[
\left\Vert P_{a}\right\Vert _{\mathcal{U}_{k}\left(  E;F\right)  }\leq
\beta_{1}\left\Vert P\right\Vert _{\mathcal{U}_{k+1}\left(  E;F\right)
}\left\Vert a\right\Vert .
\]

$\left(  ii\right)  $ For each $P\in\mathcal{U}_{k}\left(  E;F\right)  $ and
$\gamma\in E^{\ast}$, $\gamma P$ belongs to $\mathcal{U}_{k+1}\left(
E;F\right)  $ and
\[
\left\Vert \gamma P\right\Vert _{\mathcal{U}_{k+1}\left(  E;F\right)  }%
\leq\beta_{2}\left\Vert \gamma\right\Vert \left\Vert P\right\Vert
_{\mathcal{U}_{k}\left(  E;F\right)  }.
\]

\end{definition}

\section{The first multilinear and polynomial approaches to summability}

In 1989, R. Alencar and M.C. Matos \cite{am} explored the following concept of
absolutely summing multilinear operators, which was essentially introduced by Pietsch:

\begin{definition}
\label{pro}Let $p,p_{1},\ldots,p_{n}\in(0,\infty),$ with $\frac{1}{p}\leq
\frac{1}{p_{1}}+\cdots+\frac{1}{p_{n}}.$ A mapping $T\in\mathcal{L}%
(E_{1},\ldots,E_{n};F)$ is absolutely $(p;p_{1},\ldots,p_{n})$%
-summing\textbf{\ (}or\textbf{\ }$(p;p_{1},\ldots,p_{n})$-summing\textbf{)} if
there exists a $C\geq0$ such that%
\begin{equation}
\left(  \overset{m}{\underset{i=1}{\sum}}\left\Vert T(x_{i}^{(1)},\ldots
,x_{i}^{(n)})\right\Vert ^{p}\right)  ^{\frac{1}{p}}\leq C\overset
{n}{\underset{k=1}{\prod}}\left\Vert \left(  x_{j}^{(k)}\right)  _{j=1}%
^{m}\right\Vert _{w,p_{k}} \label{llo}%
\end{equation}

for every $m\in\mathbb{N}$ and $x_{i}^{(k)}\in E_{k},$ with $\left(
i,k\right)  \in\left\{  1,\ldots,m\right\}  \times\left\{  1,\ldots,n\right\}
$. Analogously an $n$-homogeneous polynomial $P\in\mathcal{P}(^{n}E;F)$
is\textbf{\ }absolutely\textbf{\ }$(p;q)$-summing if there exists a constant
$C\geq0$ such that%
\[
\left(  \sum\limits_{j=1}^{m}\left\Vert P\left(  x_{j}\right)  \right\Vert
^{p}\right)  ^{\frac{1}{p}}\leq C\left\Vert \left(  x_{j}\right)  _{j=1}%
^{m}\right\Vert _{w,q}^{n}%
\]
for all $m\in\mathbb{N}$ and $x_{j}\in E,$ with $j=1,\ldots,m.$
\end{definition}

The space of all $n$-linear operators satisfying (\ref{llo}) will be denoted
by $\mathcal{L}_{as(p;p_{1},\ldots,p_{n})}(E_{1},\ldots,E_{n};F).$ When
$p_{1}=\cdots=p_{n}=q,$ we simply write $\mathcal{L}_{as(p;q)}(E_{1}
,\ldots,E_{n};F)$. For $n=1$ we use the classical notation $\Pi_{(p;q)}$
instead of $\mathcal{L}_{as(p;q)}.$ For polynomials we write $\mathcal{P}
_{as(p;q)}(^{n}E;F).$

For other approaches we mention \cite{port, Choi, df, Dimant} and references
therein. The successful notion of multiple summing multilinear operators will
be mentioned in the Section \ref{MMMA}.

In the case of mixing summing operators, the multilinear/polynomial theory was
investigated by C.A. Soares in his Ph.D. dissertation \cite{Carlos Alberto}.
However, the definition considered in \cite{Carlos Alberto} is an extension of
Definition \ref{pro} and, as it happens to the concept of absolutely summing
multilinear operators, it inherits its weaknesses.

\begin{definition}
Let $0<q\leq s\leq\infty$ and $0<p_{1},\ldots,p_{n}\leq\infty.$ An $n$-linear
operator $T\in\mathcal{L}(E_{1},\ldots,E_{n};F)$ is\textbf{\ }$(s,q;p_{1}%
,\ldots,p_{n})$-mixing summing if there exists a constant $\sigma\geq0$ such
that
\begin{equation}
\left\Vert \left(  T(x_{j}^{(1)},\ldots,x_{j}^{(n)})\right)  _{j=1}%
^{m}\right\Vert _{mx(s,q)}\leq\sigma\prod_{k=1}^{n}\left\Vert (x_{j}%
^{(k)})_{j=1}^{m}\right\Vert _{w,p_{k}} \label{xxc}%
\end{equation}
for every $m\in\mathbb{N}$ $,$ $x_{1}^{(1)},\ldots,x_{m}^{(1)}\in E_{1}%
,\ldots,x_{1}^{(n)},\ldots,x_{m}^{(n)}\in E_{n}.$ Analogously $P\in
\mathcal{P}(^{n}E;F)$ is\textbf{\ }mixing\textbf{\ }$(s,q;p)$-summing if there
exists a constant $C\geq0$ such that
\[
\left\Vert \left(  P(x_{j})\right)  _{j=1}^{m}\right\Vert _{mx(s,q)}\leq
C\left\Vert \left(  x_{j}\right)  _{j=1}^{m}\right\Vert _{w,p}^{n}%
\]
for all $m\in\mathbb{N}$ and $x_{j}\in E,$ with $j=1,\ldots,m.$
\end{definition}

If $p_{1}=\cdots=p_{n}=p,$ the operator $T$ is said $(s,q;p)$-mixing summing.

The following multilinear generalization of $(p;q;r)$-summing operators was
recently introduced by D. Achour \cite{achour}:

\begin{definition}
\label{ach}Let $0<p,q_{1},\ldots,q_{n}<\infty$ and $0<r\leq\infty$ with
\[
\frac{1}{p}\leq\frac{1}{q_{1}}+\cdots+\frac{1}{q_{n}}+\frac{1}{r}.
\]
An $n$-linear map $T$ $\in\mathcal{L}{(E_{1},\ldots,E_{n};F)}$ is absolutely
$(p;q_{1},\ldots,q_{n};r)$-summing if there is a $C\geq0$ so that
\begin{equation}
\left(  \sum\limits_{j=1}^{m}\left\vert \varphi_{j}\left(  T\left(
x_{j}^{(1)},\ldots,x_{j}^{(n)}\right)  \right)  \right\vert ^{p}\right)
^{\frac{1}{p}}\leq C\left\Vert \left(  \varphi_{j}\right)  _{j=1}
^{m}\right\Vert _{w,r}\prod_{i=1}^{n}\left\Vert \left(  x_{j}^{(i)}\right)
_{j=1}^{m}\right\Vert _{w,q_{i}} \label{des}%
\end{equation}
for all $m\in\mathbb{N}$, $\varphi_{j}\in F^{\ast}$ and $x_{j}^{(i)}\in
E_{i},$ with $\left(  i,j\right)  \in\{1,\ldots,n\}\times\{1,\ldots,m\}.$
Analogously an $n$-homogeneous polynomial $P\in\mathcal{P}(^{n}E;F)$
is\textbf{\ }absolutely\textbf{\ }$(p;q;r)$-summing if there exists a constant
$C\geq0$ such that
\[
\left(  \sum\limits_{j=1}^{m}\left\vert \varphi_{j}\left(  P\left(
x_{j}\right)  \right)  \right\vert ^{p}\right)  ^{\frac{1}{p}}\leq C\left\Vert
\left(  \varphi_{j}\right)  _{j=1}^{m}\right\Vert _{w,r}\left\Vert \left(
x_{j}\right)  _{j=1}^{m}\right\Vert _{w,q}^{n}%
\]
for all $m\in\mathbb{N}$, $\varphi_{j}\in F^{\ast}$ and $x_{j}\in E,$ with
$j=1,\ldots,m.$
\end{definition}

We denote{\ the space of all absolutely }$(p;q_{1},\ldots,q_{n};r)$-summing
$n$-linear operators by
\[
\mathcal{L}_{as(p;q_{1},\ldots,q_{n};r)}\left(  {E_{1},\ldots,E_{n};F}\right)
.
\]
When $q_{1}=\cdots=q_{n}=q$ we just write $\mathcal{L}_{as(p;q;r)}\left(
{E_{1},\ldots,E_{n};F}\right)  $. When $r=\infty$ we recover the notion of
absolutely $\left(  p;q_{1}, \ldots,q_{n}\right)  $-summing multilinear
mappings $\mathcal{L}_{as(p;q_{1},\ldots,q_{n})}$ due to Alencar and Matos
\cite{am}. More precisely,
\begin{equation}
\mathcal{L}_{as(p;q_{1},\ldots,q_{n};\infty)}=\mathcal{L}_{as(p;q_{1}
,\ldots,q_{n})}. \label{s32}%
\end{equation}

If $\frac{1}{p}>\frac{1}{q_{1}}+\cdots+\frac{1}{q_{n}}+\frac{1}{r}$ and $T$ is
absolutely $(p;q_{1},\ldots,q_{n};r)$-summing, then $T=0$. \ It is not
difficult to prove that
\begin{equation}
\mathcal{L}_{as\left(  p;q_{1},\ldots,q_{n}\right)  }\left(  E_{1}
,\ldots,E_{n};F\right)  \subset\mathcal{L}_{as\left(  p;q_{1},\ldots
,q_{n};r\right)  }\left(  E_{1},\ldots,E_{n};F\right)  \label{mov}%
\end{equation}
for all Banach spaces $E_{1},\ldots,E_{n},F$ and $r>0$.

\subsection{The lack of coherence and compatibility}

The class of absolutely $\left(  p;q\right)  $-summing $n$-homogeneous
polynomials will be denoted by $\mathcal{P}_{as(p;q)}^{n}.$ As before, the
space of all $n$-homogeneous polynomials $P:E\rightarrow F$ in $\mathcal{P}%
_{as(p;q)}^{n}$ is represented by $\mathcal{P}_{as(p;q)}\left(  ^{n}%
E;F\right)  .$ The notions of absolutely $\left(  p;q;r\right)  $-summing
polynomials and mixing summing polynomials are denoted in a similar way.

It can be easily seen that $\left(  \mathcal{P}_{as(p;q)}^{n}\right)
_{n=1}^{\infty}$ in general fails to be coherent and compatible with
$\Pi_{as(p;q)}$. In fact for any positive integer $n\geq2$ and any real number
$1\leq p\leq2$ we know that
\[
\mathcal{P}_{as(1;1)}\left(  ^{n}\ell_{p};F\right)  =\mathcal{P}\left(
^{n}\ell_{p};F\right)
\]
for all Banach spaces $F$. This result is an obvious deviation from the spirit
of the linear ideal of absolutely summing operators since
\[
\Pi_{as(1;1)}\left(  \ell_{p};F\right)  =\mathcal{L}\left(  \ell_{p};F\right)
\]
if and only if $p=1$ and $F$ is a Hilbert space (see \cite{LP}). This
situation also proves that $\left(  \mathcal{P}_{as(1;1)}^{n}\right)
_{n=1}^{\infty}$ is not coherent or compatible with $\Pi_{as(1;1)}.$ We also
know that $\left(  \mathcal{P}_{as(p;q)}^{n}\right)  _{n=1}^{\infty}$ in
general is not a (global) holomorphy type.

Since $\mathcal{P}_{as\left(  p;q;\infty\right)  }^{n}=\mathcal{P}_{as\left(
p;q\right)  }^{n}$ and $\mathcal{P}_{mxs\left(  \infty;p\right)  }
^{n}=\mathcal{P}_{as\left(  p;p\right)  }^{n}$ these deficiencies of $\left(
\mathcal{P}_{as(1;1)}^{n}\right)  _{n=1}^{\infty}$are inherited by the
polynomial analogues of the concepts of Subsections \ref{xxzzz} and
\ref{Sub2}. These deficiencies shall be fixed by the alternative concepts
introduced in the next sections.

\section{Multiple summing multilinear operators: the ``nice
prototype''\label{MMMA}}

Multiple $(p;q)$-summing multilinear were introduced in 2003 \cite{Matos, pv}.
The origins of this notion date back to the 1930's with Littlewood's $4/3$
inequality \cite{lit} which asserts that
\[
\left(  \sum\limits_{i,j=1}^{N}\left\vert T(e_{i},e_{j})\right\vert ^{\frac
{4}{3}}\right)  ^{\frac{3}{4}}\leq\sqrt{2}\left\Vert T\right\Vert
\]
for every bilinear form $T:\ell_{\infty}^{N}\times\ell_{\infty}^{N}
\rightarrow\mathbb{K}$ and every positive integer $N.$ In 1931 H.F.
Bohnenblust and E. Hille \cite{bh} provided a deep generalization of this
result to multilinear mappings: for every positive integer $n$ there is a
$C_{n}>0$ so that
\[
\left(  \sum\limits_{i_{1},\ldots,i_{n}=1}^{N}\left\vert T(e_{i_{^{1}}}
,\ldots,e_{i_{n}})\right\vert ^{\frac{2n}{n+1}}\right)  ^{\frac{n+1}{2n}}\leq
C_{n}\left\Vert T\right\Vert
\]
for every $n$-linear mapping $T:\ell_{\infty}^{N}\times\cdots\times
\ell_{\infty}^{N}\rightarrow\mathbb{C}$ and every positive integer $N$. This
result has important applications in operator theory in Banach spaces,
harmonic analysis, complex analysis and analytic number theory. For recent
advances related to the Bohnenblust-Hille inequality we refer to \cite{an,
df,DPell, Mun, Mu2}.

In his Ph.D. dissertation, D. P\'{e}rez-Garc\'{\i}a \cite{Da} remarked that
the Bohnenblust-Hille inequality can be viewed as a result of the theory of
multiple summing operators.

\begin{theorem}
[Bohnenblust-Hille]\label{ytr}If $E_{1},\ldots,E_{n}$ are Banach spaces and
$T\in\mathcal{L}(E_{1},\ldots,E_{n};\mathbb{K}),$ then there exists a constant
$C_{n}\geq0$ such that%
\begin{equation}
\left(  \sum_{j_{1},\ldots,j_{n}=1}^{N}\left\vert T(x_{j_{1}}^{(1)}%
,\ldots,x_{j_{n}}^{(n)})\right\vert ^{\frac{2n}{n+1}}\right)  ^{\frac{n+1}%
{2n}}\leq C_{n}\prod_{k=1}^{n}\left\Vert (x_{j}^{(k)})_{j=1}^{N}\right\Vert
_{w,1} \label{juo}%
\end{equation}
for every positive integer $N$ and $x_{j}^{(k)}\in E_{k}$, $k=1,\ldots,n$ and
$j=1,\ldots,N.$
\end{theorem}

The inequality above can be regarded as a result in the theory of multiple
summing multilinear operators. Recall that for $1\leq q_{1},\ldots,q_{n}\leq
p<\infty,$ an $n$-linear operator $T:E_{1}\times\cdots\times E_{n}\rightarrow
F$ is multiple\emph{\ }$(p;q_{1},\ldots,q_{n})$-summing ($T\in\mathcal{L}%
_{mas(p;q_{1},\ldots,q_{n})}(E_{1},\ldots,E_{n};F)$) if there exists $C>0$
such that
\begin{equation}
\left(  \sum_{j_{1},\ldots,j_{n}=1}^{\infty}\Vert T(x_{j_{1}}^{(1)}%
,\ldots,x_{j_{n}}^{(n)})\Vert^{p}\right)  ^{1/p}\leq C\prod\limits_{k=1}%
^{n}\Vert(x_{j}^{(k)})_{j=1}^{\infty}\Vert_{w,q_{k}}\text{ } \label{jup2}%
\end{equation}
for every $(x_{j}^{(k)})_{j=1}^{\infty}\in\ell_{q_{k}}^{w}(E_{k})$,
$k=1,\ldots,n$.

The infimum of all $C$'s satisfying (\ref{jup2}), denoted by $\left\Vert
T\right\Vert _{(r;r_{1},\ldots,r_{n})},$ defines a complete norm if $r\geq1$
($r$-norm, if $r\in(0,1)$) in $\mathcal{L}_{mas(r;r_{1},\ldots,r_{n})}%
(E_{1},\ldots,E_{n};F).$ If $r_{1}=\cdots=r_{n}=s$ we just write $(r;s),$ and
when $r=s$ we replace $\left(  r;r\right)  $ by $r$. For $n=1$ this concept
also coincides with the classical notion of absolutely summing linear
operators and, for this reason, we keep the usual notation $\pi_{(r;s)}\left(
T\right)  $ instead of $\left\Vert T\right\Vert _{(r;s)}$ for the norm of $T.$
The essence of the notion of multiple summing multilinear operators, for
bilinear operators, can also be traced back to \cite{Ram}. For recent results
in the theory of multiple summing operators we refer to \cite{BBJP, se,
davidstudia, PopJM} and references therein.

\section{Multiple $\left(  p;q_{1},\ldots,q_{n};r\right)  $-summing
multilinear operators\label{ss33}}

In this section we introduce the notion of multiple $\left(  p;q_{1}%
,\ldots,q_{n};r\right)  $-summing multilinear operators and, as we shall see
in the next sections, the polynomial version of this concept is coherent and
compatible with the (linear) operator ideal of $(p;q;r)$-summing operators.

\begin{definition}
Let $m\in\mathbb{N},p,r,q_{1},\ldots,q_{n}\geq1$ and $E_{1},\ldots,E_{n},F$ be
Banach spaces. A continuous multilinear operator $T:E_{1}\times\cdots\times
E_{n}\rightarrow F$ is multiple $\left(  p;q_{1},\ldots,q_{n};r\right)
$-summing when%
\[
\left(  \varphi_{j_{1}\ldots j_{n}}\left(  T\left(  x_{j_{1}}^{\left(
1\right)  },\ldots,x_{j_{n}}^{\left(  n\right)  }\right)  \right)  \right)
_{j_{1},\ldots,j_{n}\in\mathbb{N}}\in\ell_{p}\left(  \mathbb{N}^{n}\right)
\]
whenever $\left(  x_{j}^{\left(  i\right)  }\right)  _{j=1}^{\infty}\in
\ell_{q_{i}}^{w}\left(  E_{i}\right)  ,i=1,\ldots,n$ and $\left(
\varphi_{j_{1}\ldots j_{n}}\right)  _{j_{1},\ldots,j_{n}\in\mathbb{N}}\in
\ell_{r}^{w}\left(  F^{\ast},\mathbb{N}^{n}\right)  .$
\end{definition}

Sometimes we shall simply write $j\in\mathbb{N}^{n}$ to denote $j=(j_{1}%
,\ldots,j_{n})\in\mathbb{N}^{n}.$ The vector space formed by the multiple
$\left(  p;q_{1},\ldots,q_{n};r\right)  $-summing multilinear operators from
$E_{1}\times\cdots\times E_{n}$ to $F$ shall be represented by $\mathcal{L}%
_{mas\left(  p;q_{1},\ldots,q_{n};r\right)  }\left(  E_{1},\ldots
,E_{n};F\right)  $). When $q_{1}=\cdots=q_{n}=q$, we simply write
$\mathcal{L}_{mas\left(  p;q;r\right)  }\left(  E_{1},\ldots,E_{n};F\right)  $.

As it happens in other similar classes, the class $\mathcal{L}_{mas\left(
p;q_{1},\ldots,q_{n};r\right)  }\left(  E_{1},\ldots,E_{n};F\right)  $ has a
characterization by means of inequalities:

\begin{theorem}
\label{1.2}The following assertions are equivalent for $T\in\mathcal{L}\left(
E_{1},\ldots,E_{n};F\right)  $:

\begin{itemize}
\item[(i)] $T\in\mathcal{L}_{mas\left(  p;q_{1},\ldots,q_{n};r\right)
}\left(  E_{1},\ldots,E_{n};F\right)  ;$

\item[(ii)] There is a $C\geq0$ such that
\begin{align}
&  \left(  \sum_{j_{1},\ldots,j_{n}=1}^{\infty}\left\vert \varphi_{j_{1}\ldots
j_{n}}\left(  T\left(  x_{j_{1}}^{\left(  1\right)  },\ldots,x_{j_{n}%
}^{\left(  n\right)  }\right)  \right)  \right\vert ^{p}\right)  ^{\frac{1}%
{p}}\label{11:06}\\
&  \leq C\left\Vert \left(  \varphi_{j_{1}\ldots j_{n}}\right)  _{j_{1}%
,...,j_{n}\in\mathbb{N}}\right\Vert _{w,r}\prod_{i=1}^{n}\left\Vert \left(
x_{j}^{\left(  i\right)  }\right)  _{j=1}^{\infty}\right\Vert _{w,q_{i}%
}\nonumber
\end{align}
whenever $\left(  x_{j}^{\left(  i\right)  }\right)  _{j=1}^{\infty}\in
\ell_{q_{i}}^{w}\left(  E_{i}\right)  ,i=1,\ldots,n$ and $\left(
\varphi_{j_{1}\ldots j_{n}}\right)  _{j\in\mathbb{N}^{n}}\in\ell_{r}%
^{w}\left(  F^{\ast},\mathbb{N}^{n}\right)  ;$

\item[(iii)] There is a $C\geq0$ such that
\begin{align*}
&  \left(  \sum_{j_{1},\ldots,j_{n}=1}^{m}\left\vert \varphi_{j_{1}\ldots
j_{n}}\left(  T\left(  x_{j_{1}}^{\left(  1\right)  },\ldots,x_{j_{n}%
}^{\left(  n\right)  }\right)  \right)  \right\vert ^{p}\right)  ^{\frac{1}%
{p}}\\
&  \leq C\left\Vert \left(  \varphi_{j_{1}\ldots j_{n}}\right)  _{j_{1}%
,...,j_{n}\in\mathbb{N}_{m}}\right\Vert _{w,r}\prod_{i=1}^{n}\left\Vert
\left(  x_{j}^{\left(  i\right)  }\right)  _{j=1}^{m}\right\Vert _{w,q_{i}}%
\end{align*}
for all $m\in\mathbb{N},$ $x_{1}^{\left(  i\right)  },\ldots,x_{m}^{\left(
i\right)  }\in E_{i},i=1,\ldots,n$ and $\left(  \varphi_{j_{1}\ldots j_{n}%
}\right)  _{j\in\mathbb{N}_{m}^{n}}\in\ell_{r}^{w}\left(  F^{\ast}%
,\mathbb{N}_{m}^{n}\right)  .$
\end{itemize}

The infimum of all $C$ satisfying (\ref{11:06}) defines a norm in
$\mathcal{L}_{mas\left(  p;q_{1},\ldots,q_{n};r\right)  }\left(  E_{1}%
,\ldots,E_{n};F\right)  .$
\end{theorem}

Similarly to (\ref{mov}) it can also be proved that
\begin{equation}
\mathcal{L}_{mas\left(  p;q_{1},\ldots,q_{n}\right)  }\subset\mathcal{L}
_{mas\left(  p;q_{1},\ldots,q_{n};r\right)  } \label{inc222}%
\end{equation}
for all $r>0.$ From Theorem \ref{1.2} we can conclude that if
\[
\frac{1}{p}>\frac{1}{q_{i}}+\frac{1}{r}%
\]
for some $i$, then $\mathcal{L}_{mas\left(  p;q_{1},\ldots,q_{n};r\right)
}\left(  E_{1},\ldots,E_{n};F\right)  =\left\{  0\right\}  $. In fact, we
first prove that if $T\in\mathcal{L}_{mas\left(  p;q_{1},\ldots,q_{n}
;r\right)  }\left(  E_{1},\ldots,E_{n};F\right)  ,$ then, for any $a\in E_{1}
$, the map
\begin{equation}
T_{a}:E_{2}\times\cdots\times E_{n}\longrightarrow F:T_{a}\left(  x_{2}
,\ldots,x_{n}\right)  =T\left(  a,x_{2},\ldots,x_{n}\right)  \label{wsq}%
\end{equation}
is multiple $\left(  p;q_{2},\ldots,q_{n};r\right)  $-summing and
\begin{equation}
\left\Vert T\right\Vert _{mas\left(  p;q_{2},\ldots,q_{n};r\right)  }
\leq\left\Vert a\right\Vert \left\Vert T\right\Vert _{mas\left(
p;q_{1},\ldots,q_{n};r\right)  }. \label{wsa}%
\end{equation}
So, if $\frac{1}{p}>\frac{1}{q_{i}}+\frac{1}{r}$ for some $i$, then
$\mathcal{L}_{mas\left(  p;q_{1},\ldots,q_{n};r\right)  }\left(  E_{1}%
,\ldots,E_{n};F\right)  =\left\{  0\right\}  $. In fact, suppose that
$\frac{1}{p}>\frac{1}{q_{1}}+\frac{1}{r}.$ So, using \ (\ref{wsq}), we know
that if $T\in\mathcal{L}_{mas\left(  p;q_{1},\ldots,q_{n};r\right)  }\left(
E_{1},\ldots,E_{n};F\right)  $ then $T_{a_{2},\ldots,a_{n}}\in\mathcal{L}%
_{as\left(  p;q_{1};r\right)  }\left(  E_{1};F\right)  $ for all $a_{2}\in
E_{2},\ldots,a_{n}\in E_{n}$. It follows that $T_{a_{2},\ldots,a_{n}}=0$ and
hence $T=0.$ So, in order to avoid trivialities we shall suppose $\frac{1}%
{p}\leq\frac{1}{q_{i}}+\frac{1}{r}$ for all $i.$

\subsection{Coherence and compatibility \label{ss44}}

Standard calculations show that
\[
\left(  \mathcal{L}_{mas\left(  p;q_{1},\ldots,q_{n};r\right)  },\left\Vert
\cdot\right\Vert _{mas\left(  p;q_{1},\ldots,q_{n};r\right)  }\right)
\]
is a Banach multi-ideal. If $\mathcal{M}$ is a (quasi-) normed ideal of
multilinear mappings, the class
\[
\mathcal{P}_{\mathcal{M}}=\left\{  P\in\mathcal{P}^{n};\check{P}\in
\mathcal{M},n\in\mathbb{N}\right\}  \text{,}%
\]
with $\left\Vert P\right\Vert _{\mathcal{P}_{\mathcal{M}}}:=\left\Vert
\check{P}\right\Vert _{\mathcal{M}},$ is a (quasi-) normed ideal of
polynomials, called polynomial ideal generated by $\mathcal{M}$. If
$\mathcal{M}$ is (quasi-) Banach, then $\mathcal{P}_{\mathcal{M}}$ is (quasi-)
Banach (see \cite[p. 46]{BBJMs}).

Thus, the class
\[
\mathcal{P}_{mas\left(  p;q;r\right)  }^{n}=\left\{  P\in\mathcal{P}%
^{n};\check{P}\in\mathcal{L}_{mas\left(  p;q;r\right)  }^{n}\right\}  ,
\]
with
\[
\left\Vert P\right\Vert _{\mathcal{P}_{mas\left(  p;q;r\right)  }^{n}%
}:=\left\Vert \check{P}\right\Vert _{mas\left(  p;q;r\right)  },
\]
ia a Banach polynomial ideal$.$

\begin{theorem}
$\left(  \mathcal{P}_{mas\left(  p;q;r\right)  }^{n},\left\Vert .\right\Vert
_{\mathcal{P}_{mas\left(  p;q;r\right)  }^{n}}\right)  _{n=1}^{\infty}$ is
coherent and, for each fixed $n$, compatible with $\mathcal{L}_{mas\left(
p;q;r\right)  }$.
\end{theorem}

\begin{proof}
If $P\in\mathcal{P}_{mas\left(  p;q;r\right)  }^{n}\left(  ^{n}E;F\right)  $
and $a\in E$, then $\check{P}\in\mathcal{L}_{mas\left(  p;q;r\right)  }%
^{n}\left(  ^{n}E;F\right)  $ and, from (\ref{wsq}) and (\ref{wsa}),
$\check{P}_{a}\in\mathcal{L}_{mas\left(  p;q;r\right)  }^{n-1}\left(
^{n-1}E;F\right)  .$ Hence $P_{a}\in\mathcal{P}_{mas\left(  p;q;r\right)
}^{n-1}\left(  ^{n-1}E;F\right)  $ with
\[
\left\Vert P_{a}\right\Vert _{\mathcal{P}_{mas\left(  p;q;r\right)  }^{n-1}%
}\leq\left\Vert a\right\Vert \left\Vert P\right\Vert _{\mathcal{P}_{mas\left(
p;q;r\right)  }^{n}}.
\]
Let $\gamma\in E^{\ast}.$ Note that%
\[
\left(  \gamma P\right)  ^{\vee}\left(  x_{1},\ldots,x_{n+1}\right)  =\frac
{1}{n+1}\sum_{k=1}^{n+1}\gamma\left(  x_{k}\right)  \check{P}\left(
x_{1},\overset{\left[  k\right]  }{\ldots},x_{n+1}\right)  ,
\]
where $\overset{\left[  k\right]  }{\ldots}$ means that the $k$-th coordinate
is missing.

Let $m\in\mathbb{N}$, $x_{j}^{(k)}\in E$, with $j=1,\ldots.,m$ and
$k=1,\ldots,n+1;$ let $\varphi_{j_{1}\ldots j_{n+1}}\in F^{\ast}$ with
$j_{1},\ldots,j_{n+1}=1,\ldots.,m.$ Using the triangle inequality we have%
\begin{align*}
&  \left(  \sum_{j_{1},\ldots,j_{n+1}=1}^{m}\left\vert \varphi_{j_{1}\ldots
j_{n+1}}\left(  \left(  \gamma P\right)  ^{\vee}\left(  x_{j_{1}}^{\left(
1\right)  },\ldots,x_{j_{n+1}}^{\left(  n+1\right)  }\right)  \right)
\right\vert ^{p}\right)  ^{\frac{1}{p}}\\
&  =\left(  \sum_{j_{1},\ldots,j_{n+1}=1}^{m}\left\vert \varphi_{j_{1}\ldots
j_{n+1}}\left(  \frac{1}{n+1}\sum_{k=1}^{n+1}\gamma\left(  x_{j_{k}}^{\left(
k\right)  }\right)  \check{P}\left(  x_{j_{1}}^{\left(  1\right)  }%
,\overset{\left[  k\right]  }{\ldots},x_{j_{n+1}}^{\left(  n+1\right)
}\right)  \right)  \right\vert ^{p}\right)  ^{\frac{1}{p}}\\
&  =\frac{1}{n+1}\left(  \sum_{j_{1},\ldots,j_{n+1}=1}^{m}\left\vert
\varphi_{j_{1}\ldots j_{n+1}}\left(  \sum_{k=1}^{n+1}\gamma\left(  x_{j_{k}%
}^{\left(  k\right)  }\right)  \check{P}\left(  x_{j_{1}}^{\left(  1\right)
},\overset{\left[  k\right]  }{\ldots},x_{j_{n+1}}^{\left(  n+1\right)
}\right)  \right)  \right\vert ^{p}\right)  ^{\frac{1}{p}}\\
&  =\frac{1}{n+1}\left(  \sum_{j_{1},\ldots,j_{n+1}=1}^{m}\left\vert
\sum_{k=1}^{n+1}\varphi_{j_{1}\ldots j_{n+1}}\left(  \gamma\left(  x_{j_{k}%
}^{\left(  k\right)  }\right)  \check{P}\left(  x_{j_{1}}^{\left(  1\right)
},\overset{\left[  k\right]  }{\ldots},x_{j_{n+1}}^{\left(  n+1\right)
}\right)  \right)  \right\vert ^{p}\right)  ^{\frac{1}{p}}\\
&  \leq\frac{1}{n+1}\left(  \sum_{j_{1},\ldots,j_{n+1}=1}^{m}\left(
\sum_{k=1}^{n+1}\left\vert \varphi_{j_{1}\ldots j_{n+1}}\left(  \gamma\left(
x_{j_{k}}^{\left(  k\right)  }\right)  \check{P}\left(  x_{j_{1}}^{\left(
1\right)  },\overset{\left[  k\right]  }{\ldots},x_{j_{n+1}}^{\left(
n+1\right)  }\right)  \right)  \right\vert \right)  ^{p}\right)  ^{\frac{1}%
{p}}\\
&  =\frac{1}{n+1}\left\Vert \left(  \sum_{k=1}^{n+1}\left\vert \varphi
_{j_{1}\ldots j_{n+1}}\left(  \gamma\left(  x_{j_{k}}^{\left(  k\right)
}\right)  \check{P}\left(  x_{j_{1}}^{\left(  1\right)  },\overset{\left[
k\right]  }{\ldots},x_{j_{n+1}}^{\left(  n+1\right)  }\right)  \right)
\right\vert \right)  _{j_{1},\ldots,j_{n+1}=1}^{m}\right\Vert _{p}\\
&  =(\ast).
\end{align*}
Thus, from the Minkowski inequality we have
\begin{align}
(\ast)  &  =\nonumber\\
&  =\frac{1}{n+1}\left\Vert \sum_{k=1}^{n+1}\left(  \left\vert \varphi
_{j_{1}\ldots j_{n+1}}\left(  \gamma\left(  x_{j_{k}}^{\left(  k\right)
}\right)  \check{P}\left(  x_{j_{1}}^{\left(  1\right)  },\overset{\left[
k\right]  }{\ldots},x_{j_{n+1}}^{\left(  n+1\right)  }\right)  \right)
\right\vert \right)  _{j_{1},\ldots,j_{n+1}=1}^{m}\right\Vert _{p}\\
&  \leq\frac{1}{n+1}\sum_{k=1}^{n+1}\left\Vert \left(  \left\vert
\varphi_{j_{1}\ldots j_{n+1}}\left(  \gamma\left(  x_{j_{k}}^{\left(
k\right)  }\right)  \check{P}\left(  x_{j_{1}}^{\left(  1\right)  }%
,\overset{\left[  k\right]  }{\ldots},x_{j_{n+1}}^{\left(  n+1\right)
}\right)  \right)  \right\vert \right)  _{j_{1},\ldots,j_{n+1}=1}%
^{m}\right\Vert _{p}\nonumber\\
&  =\frac{1}{n+1}\sum_{k=1}^{n+1}\left(  \sum_{j_{1},\ldots,j_{n+1}=1}%
^{m}\left\vert \varphi_{j_{1}\ldots j_{n+1}}\left(  \gamma\left(  x_{j_{k}%
}^{\left(  k\right)  }\right)  \check{P}\left(  x_{j_{1}}^{\left(  1\right)
},\overset{\left[  k\right]  }{\ldots},x_{j_{n+1}}^{\left(  n+1\right)
}\right)  \right)  \right\vert ^{p}\right)  ^{\frac{1}{p}}\nonumber\\
&  =\frac{1}{n+1}\left[  \left(  \sum_{j_{1},\ldots,j_{n+1}=1}^{m}\left\vert
\varphi_{j_{1}\ldots j_{n+1}}\left(  \check{P}\left(  \gamma\left(  x_{j_{1}%
}^{\left(  1\right)  }\right)  x_{j_{2}}^{\left(  2\right)  },\ldots
,x_{j_{n+1}}^{\left(  n+1\right)  }\right)  \right)  \right\vert ^{p}\right)
^{\frac{1}{p}}+\cdots\right. \nonumber\\
&  \left.  \cdots+\left(  \sum_{j_{1},\ldots,j_{n+1}=1}^{m}\left\vert
\varphi_{j_{1}\ldots j_{n+1}}\left(  \check{P}\left(  \gamma\left(
x_{j_{n+1}}^{\left(  n+1\right)  }\right)  x_{j_{1}}^{\left(  1\right)
},\ldots,x_{j_{n}}^{\left(  n\right)  }\right)  \right)  \right\vert
^{p}\right)  ^{\frac{1}{p}}\right]  .\nonumber
\end{align}

Hence
\begin{align}
&  \left(  \sum_{j_{1},\ldots,j_{n+1}=1}^{m}\left\vert \varphi_{j_{1}\ldots
j_{n+1}}\left(  \left(  \gamma P\right)  ^{\vee}\left(  x_{j_{1}}^{\left(
1\right)  },\ldots,x_{j_{n+1}}^{\left(  n+1\right)  }\right)  \right)
\right\vert ^{p}\right)  ^{\frac{1}{p}}\label{estta}\\
&  \leq\frac{1}{n+1}\left[  \left(  \sum_{j_{1},\ldots,j_{n+1}=1}%
^{m}\left\vert \varphi_{j_{1}\ldots j_{n+1}}\left(  \check{P}\left(
\gamma\left(  x_{j_{1}}^{\left(  1\right)  }\right)  x_{j_{2}}^{\left(
2\right)  },\ldots,x_{j_{n+1}}^{\left(  n+1\right)  }\right)  \right)
\right\vert ^{p}\right)  ^{\frac{1}{p}}+\cdots\right. \nonumber\\
&  \left.  \cdots+\left(  \sum_{j_{1},\ldots,j_{n+1}=1}^{m}\left\vert
\varphi_{j_{1}\ldots j_{n+1}}\left(  \check{P}\left(  \gamma\left(
x_{j_{n+1}}^{\left(  n+1\right)  }\right)  x_{j_{1}}^{\left(  1\right)
},\ldots,x_{j_{n}}^{\left(  n\right)  }\right)  \right)  \right\vert
^{p}\right)  ^{\frac{1}{p}}\right]  .\nonumber
\end{align}

Note that each one of the $n+1$ terms of (\ref{estta}) can be re-written as
\[
\left(  \sum_{j_{2}=1}^{m^{2}}\sum_{j_{3},\ldots,j_{n+1}=1}^{m}\left\vert
\widetilde{\varphi}_{j_{2}\ldots j_{n+1}}\left(  \check{P}\left(  z_{j_{2}%
}^{(2)},\ldots,z_{j_{n+1}}^{\left(  n+1\right)  }\right)  \right)  \right\vert
^{p}\right)  ^{\frac{1}{p}}%
\]
for adequate choices of $\widetilde{\varphi}_{j_{2}\ldots j_{n+1}}$ and
$z_{j_{k}}^{(k)}$, with $k=2,\ldots,n+1.$

In fact, for
\[
\left(  \sum_{j_{1},\ldots,j_{n+1}=1}^{m}\left\vert \varphi_{j_{1}\ldots
j_{n+1}}\left(  \check{P}\left(  \gamma\left(  x_{j_{1}}^{\left(  1\right)
}\right)  x_{j_{2}}^{\left(  2\right)  },\ldots,x_{j_{n+1}}^{\left(
n+1\right)  }\right)  \right)  \right\vert ^{p}\right)  ^{\frac{1}{p}},
\]
we choose%
\[
\left\{
\begin{array}
[c]{c}%
z_{j_{2}}^{\left(  2\right)  }=\gamma\left(  x_{1}^{\left(  1\right)
}\right)  x_{j_{2}}^{\left(  2\right)  }\text{ for all }j_{2}=1,\ldots.,m,\\
z_{m+j_{2}}^{\left(  2\right)  }=\gamma\left(  x_{2}^{\left(  1\right)
}\right)  x_{j_{2}}^{\left(  2\right)  }\text{ for all }j_{2}=1,\ldots.,m,\\
\vdots\\
z_{\left(  m-1\right)  m+j_{2}}^{\left(  2\right)  }=\gamma\left(
x_{m}^{\left(  1\right)  }\right)  x_{j_{2}}^{\left(  2\right)  }\text{ for
all }j_{2}=1,\ldots.,m,\\
z_{j_{i}}^{(i)}=x_{j_{i}}^{\left(  i\right)  }\text{ for all }j_{i}%
=1,\ldots,m,i=3,\ldots,n+1
\end{array}
\right.
\]
and%
\[
\left\{
\begin{array}
[c]{c}%
\widetilde{\varphi}_{j_{2},\ldots.j_{n+1}}=\varphi_{1j_{2}\ldots j_{n+1}%
}\text{ for all }j_{2}=1,\ldots.,m,\\
\widetilde{\varphi}_{m+j_{2},\ldots.j_{n+1}}=\varphi_{2j_{2}\ldots j_{n+1}%
}\text{ for all }j_{2}=1,\ldots.,m,\\
\vdots\\
\widetilde{\varphi}_{(m-1)m+j_{2},\ldots.j_{n+1}}=\varphi_{mj_{2}\ldots
j_{n+1}}\text{ for all }j_{2}=1,\ldots.,m.
\end{array}
\right.
\]
For these choices one can check that%
\begin{align*}
&  \left(  \sum_{j_{1},\ldots,j_{n+1}=1}^{m}\left\vert \varphi_{j_{1}\ldots
j_{n+1}}\left(  \check{P}\left(  \gamma\left(  x_{j_{1}}^{\left(  1\right)
}\right)  x_{j_{2}}^{\left(  2\right)  },\ldots,x_{j_{n+1}}^{\left(
n+1\right)  }\right)  \right)  \right\vert ^{p}\right)  ^{\frac{1}{p}}\\
&  =\left(  \sum_{j_{2}=1}^{m^{2}}\sum_{j_{3},\ldots,j_{n+1}=1}^{m}\left\vert
\widetilde{\varphi}_{j_{2}\ldots j_{n+1}}\left(  \check{P}\left(  z_{j_{2}%
}^{(2)},\ldots,z_{j_{n+1}}^{\left(  n+1\right)  }\right)  \right)  \right\vert
^{p}\right)  ^{\frac{1}{p}}
\end{align*}
and the other cases are similar. Then
\begin{align*}
&  \left(  \sum_{j_{1},\ldots,j_{n+1}=1}^{m}\left\vert \varphi_{j_{1}\ldots
j_{n+1}}\left(  \check{P}\left(  \gamma\left(  x_{j_{1}}^{\left(  1\right)
}\right)  x_{j_{2}}^{\left(  2\right)  },\ldots,x_{j_{n+1}}^{\left(
n+1\right)  }\right)  \right)  \right\vert ^{p}\right)  ^{\frac{1}{p}}\\
&  =\left(  \sum_{j_{2},\ldots,j_{n+1}=1}^{m^{2},m,\ldots,m}\left\vert
\widetilde{\varphi}_{j_{2}\ldots j_{n+1}}\left(  \check{P}\left(  z_{j_{2}%
}^{(2)},\ldots,z_{j_{n+1}}^{\left(  n+1\right)  }\right)  \right)  \right\vert
^{p}\right)  ^{\frac{1}{p}}\\
&  \leq\left\Vert \check{P}\right\Vert _{mas\left(  p;q;r\right)  }\left\Vert
\left(  \widetilde{\varphi}_{j_{2}\ldots j_{n+1}}\right)  _{j_{2}%
,\ldots,j_{n+1}}^{m^{2},m,\ldots,m}\right\Vert _{w,r}\left\Vert \left(
z_{j_{2}}^{\left(  2\right)  }\right)  _{j_{2}=1}^{m^{2}}\right\Vert
_{w,q}\prod_{i=3}^{n+1}\left\Vert \left(  z_{j_{i}}^{\left(  i\right)
}\right)  _{j_{i}=1}^{m}\right\Vert _{w,q}\\
&  =\left\Vert \check{P}\right\Vert _{mas\left(  p;q;r\right)  }\left\Vert
\left(  \varphi_{j_{1}\ldots j_{n+1}}\right)  _{j\in\mathbb{N}_{m}^{n+1}%
}\right\Vert _{w,r}\left\Vert \left(  \gamma\left(  x_{j_{1}}^{\left(
1\right)  }\right)  x_{j_{2}}^{\left(  2\right)  }\right)  _{j_{1},j_{2}%
=1}^{m}\right\Vert _{w,q}\prod_{i=3}^{n+1}\left\Vert \left(  x_{j}^{\left(
i\right)  }\right)  _{j=1}^{m}\right\Vert _{w,q}.
\end{align*}
Since%
\begin{align*}
&  \left\Vert \left(  \gamma\left(  x_{j_{1}}^{\left(  1\right)  }\right)
x_{j_{2}}^{\left(  2\right)  }\right)  _{j_{1},j_{2}=1}^{m}\right\Vert
_{w,q}\\
&  \leq\left\Vert \left(  \gamma\left(  x_{j_{1}}^{\left(  1\right)  }\right)
\right)  _{j_{1}=1}^{m}\right\Vert _{\infty}\sup_{\left\Vert \varphi
\right\Vert \leq1}\left(  \sum_{j=1}^{m}\left\vert \varphi\left(  x_{j_{2}%
}^{\left(  2\right)  }\right)  \right\vert ^{q}\right)  ^{\frac{1}{q}}\\
&  \leq\left\Vert \left(  \gamma\left(  x_{j_{1}}^{\left(  1\right)  }\right)
\right)  _{j_{1}=1}^{m}\right\Vert _{q}\left\Vert \left(  x_{j_{2}}^{\left(
2\right)  }\right)  _{j_{2}=1}^{m}\right\Vert _{w,q}\\
&  \leq\left\Vert \gamma\right\Vert \left\Vert \left(  x_{j_{1}}^{\left(
1\right)  }\right)  _{j_{1}=1}^{m}\right\Vert _{w,q}\left\Vert \left(
x_{j_{2}}^{\left(  2\right)  }\right)  _{j_{2}=1}^{m}\right\Vert _{w,q},
\end{align*}
we have%
\begin{align*}
&  \left(  \sum_{j_{1},\ldots,j_{n+1}=1}^{m}\left\vert \varphi_{j_{1}\ldots
j_{n+1}}\left(  \check{P}\left(  \gamma\left(  x_{j_{1}}^{\left(  1\right)
}\right)  x_{j_{2}}^{\left(  2\right)  },\ldots,x_{j_{n+1}}^{\left(
n+1\right)  }\right)  \right)  \right\vert ^{p}\right)  ^{\frac{1}{p}}\\
&  \leq\left\Vert \gamma\right\Vert \left\Vert \check{P}\right\Vert
_{mas\left(  p;q;r\right)  }\left\Vert \left(  \varphi_{j_{1}\ldots j_{n+1}%
}\right)  _{j\in\mathbb{N}_{m}^{n+1}}\right\Vert _{w,r}\prod_{i=1}%
^{n+1}\left\Vert \left(  x_{j}^{\left(  i\right)  }\right)  _{j=1}%
^{m}\right\Vert _{w,q}.
\end{align*}
Using the same idea for the other $n$ terms of (\ref{estta}), we obtain
\begin{align*}
&  \left(  \sum_{j_{1},\ldots,j_{n+1}=1}^{m}\left\vert \varphi_{j_{1}\ldots
j_{n+1}}\left(  \check{P}\left(  \gamma\left(  x_{j_{2}}^{\left(  2\right)
}\right)  x_{j_{1}}^{\left(  1\right)  },x_{j_{3}}^{\left(  3\right)  }%
\ldots,x_{j_{n+1}}^{\left(  n+1\right)  }\right)  \right)  \right\vert
^{p}\right)  ^{\frac{1}{p}}\\
&  \leq\left\Vert \gamma\right\Vert \left\Vert \check{P}\right\Vert
_{mas\left(  p;q;r\right)  }\left\Vert \left(  \varphi_{j_{1}\ldots j_{n+1}%
}\right)  _{j\in\mathbb{N}_{m}^{n+1}}\right\Vert _{w,r}\prod_{i=1}%
^{n+1}\left\Vert \left(  x_{j}^{\left(  i\right)  }\right)  _{j=1}%
^{m}\right\Vert _{w,q},
\end{align*}%
\[
\vdots
\]%
\begin{align*}
&  \left(  \sum_{j_{1},\ldots,j_{n+1}=1}^{m}\left\vert \varphi_{j_{1}\ldots
j_{n+1}}\left(  \check{P}\left(  \gamma\left(  x_{j_{n+1}}^{\left(
n+1\right)  }\right)  x_{j_{1}}^{\left(  1\right)  },x_{j_{2}}^{\left(
2\right)  }\ldots,x_{j_{n}}^{n}\right)  \right)  \right\vert ^{p}\right)
^{\frac{1}{p}}\\
&  \leq\left\Vert \gamma\right\Vert \left\Vert \check{P}\right\Vert
_{mas\left(  p;q;r\right)  }\left\Vert \left(  \varphi_{j_{1}\ldots j_{n+1}%
}\right)  _{j\in\mathbb{N}_{m}^{n+1}}\right\Vert _{w,r}\prod_{i=1}%
^{n+1}\left\Vert \left(  x_{j}^{\left(  i\right)  }\right)  _{j=1}%
^{m}\right\Vert _{w,q}.
\end{align*}
Therefore%
\begin{align*}
&  \left(  \sum_{j_{1},\ldots,j_{n+1}=1}^{m}\left\vert \varphi_{j_{1}\ldots
j_{n+1}}\left(  \left(  \gamma P\right)  ^{\vee}\left(  x_{j_{1}}^{\left(
1\right)  },\ldots,x_{j_{n+1}}^{\left(  n+1\right)  }\right)  \right)
\right\vert ^{p}\right)  ^{\frac{1}{p}}\\
&  \leq\frac{1}{n+1}\left[  \left\Vert \gamma\right\Vert \left\Vert \check
{P}\right\Vert _{mas\left(  p;q;r\right)  }\left\Vert \left(  \varphi
_{j_{1}\ldots j_{n+1}}\right)  _{j\in\mathbb{N}_{m}^{n+1}}\right\Vert
_{w,r}\prod_{i=1}^{n+1}\left\Vert \left(  x_{j}^{\left(  i\right)  }\right)
_{j=1}^{m}\right\Vert _{w,q}+\cdots\right. \\
&  \left.  \cdots+\left\Vert \gamma\right\Vert \left\Vert \check{P}\right\Vert
_{mas\left(  p;q;r\right)  }\left\Vert \left(  \varphi_{j_{1}\ldots j_{n+1}%
}\right)  _{j\in\mathbb{N}_{m}^{n+1}}\right\Vert _{w,r}\prod_{i=1}%
^{n+1}\left\Vert \left(  x_{j}^{\left(  i\right)  }\right)  _{j=1}%
^{m}\right\Vert _{w,q}\right] \\
&  =\left\Vert \gamma\right\Vert \left\Vert \check{P}\right\Vert _{mas\left(
p;q;r\right)  }\left\Vert \left(  \varphi_{j_{1}\ldots j_{n+1}}\right)
_{j\in\mathbb{N}_{m}^{n+1}}\right\Vert _{w,r}\prod_{i=1}^{n+1}\left\Vert
\left(  x_{j}^{\left(  i\right)  }\right)  _{j=1}^{m}\right\Vert _{w,q}.
\end{align*}
Finally we conclude that $\gamma P$ is multiple $\left(  p;q;r\right)
$-summing and
\begin{align*}
\left\Vert \gamma P\right\Vert _{\mathcal{P}_{mas\left(  p;q;r\right)  }%
^{n+1}}  &  \leq\left\Vert \gamma\right\Vert \left\Vert \check{P}\right\Vert
_{mas\left(  p;q;r\right)  }\\
&  =\left\Vert \gamma\right\Vert \left\Vert P\right\Vert _{\mathcal{P}%
_{mas\left(  p;q;r\right)  }^{n}}.
\end{align*}
The items (i) and (ii) from Definition \ref{IdeaisCompativeis} are obtained in
a similar way.
\end{proof}

\section{Multiple mixing summing operators}

In this section we introduce the notion of multiple mixing summing multilinear
operators (and polynomials) which is coherent and compatible with the
respective operator ideal. As another indicator that this is a correct
approach to nonlinear mixing summability, we prove a quotient theorem for
multilinear operators similar to the one for mixing summing linear operators.

\begin{definition}
Let $0<p_{1},\ldots,p_{n}\leq q\leq s<\infty$ . An $n$-linear operator
$A\in\mathcal{L}(E_{1},\ldots,E_{n};F)$ is multiple $(s,q;p_{1},\ldots,p_{n}%
)$-mixing summing if there exists a constant $\sigma\geq0$ such that
\begin{equation}
\left\Vert \left(  A(x_{j_{1}}^{(1)},\ldots,x_{j_{n}}^{(n)})\right)
_{j_{1},\ldots,j_{n}=1}^{m}\right\Vert _{mx(s,q)}\leq\sigma\prod_{k=1}%
^{n}\left\Vert (x_{j}^{(k)})_{j=1}^{m}\right\Vert _{w,p_{k}} \label{III}%
\end{equation}
for every $m\in\mathbb{N}$ $,$ $x_{1}^{(1)},\ldots,x_{m}^{(1)}\in E_{1}%
,\ldots,x_{1}^{(n)},\ldots,x_{m}^{(n)}\in E_{n}.$
\end{definition}

In this case we define
\[
\left\Vert A\right\Vert _{mx(s,q;p_{1},\ldots,p_{n})}=\inf\sigma.
\]
If $p_{1}=\cdots=p_{n}=p,$ we say that $A$ is multiple $(s,q;p)$-mixing
summing. The space of all multiple $(s,q;p_{1},\ldots,p_{n})$-mixing summing
is represented by $\Pi_{mx(s,q;p_{1},\ldots,p_{n})}.$

In order to avoid trivialities in the definition of multiple $(s,q;p_{1}
,\ldots,p_{n})$ mixing summing operators,\ we assume that $p_{k}\leq q$, for
all $k=1,\ldots,n.$ In fact, one can check that if $T\in\mathcal{L}(E_{1}
,\ldots,E_{n};F)$ is multiple $(s,q;p_{1},\ldots,p_{n})$ mixing summing and
$q<p_{k},$ for some $k,$ then $T=0.$

The following result, whose proof is standard and we omit, characterizes
multiple $(s,q;p_{1},\ldots,p_{n})$ mixing summing operators as those which
take adequate weakly summable sequences into adequate mixed summable sequences:

\begin{proposition}
\label{Primeira Prop}Let $0<p_{1},\ldots,p_{n}\leq q\leq s<\infty.$ An
operator $A\in\mathcal{L}(E_{1},\ldots,E_{n};F)$ is multiple $(s,q;p_{1}%
,\ldots,p_{n} )$-mixing summing if, and only if,
\[
\left(  A(x_{j_{1}}^{(1)},\ldots,x_{j_{n}}^{(n)})\right)  _{j_{1},\ldots,j_{n}
=1}^{\infty}\in\ell_{(s,q)}^{mx}\left(  F,\mathbb{N}^{n}\right)
\]
regardless of the choice of $(x_{i}^{(1)})_{i=1}^{\infty}\in\ell_{p_{1}}
^{w}(E_{1}),\ldots,$ $(x_{i}^{(n)})_{i=1}^{\infty}$ $\in\ell_{p_{n}}^{w}%
(E_{n}).$
\end{proposition}

In fact the proof of the previous proposition also shows that $A$ is multiple
$(s,q;p_{1},\ldots,p_{n})$-mixing summing if, and only if, the $n$-linear
operator
\[
\tilde{A}\left(  (x_{i}^{(1)})_{i=1}^{\infty},\ldots,(x_{i}^{(n)}%
)_{i=1}^{\infty}\right)  =\left(  A(x_{j_{1}}^{(1)},\ldots,x_{j_{n}}%
^{(n)})\right)  _{j_{1},...,j_{n}=1}^{\infty}%
\]
belongs to $\mathcal{L}(\ell_{p_{1}}^{w}(E_{1}),\ldots,\ell_{p_{n}}^{w}%
(E_{n});\ell_{(s,q)}^{mx}\left(  F,\mathbb{N}^{n}\right)  )$. Moreover
\[
\left\Vert A\right\Vert _{mx(s,q;p_{1},\ldots,p_{n})}=\left\Vert \tilde
{A}\right\Vert .
\]
The main result of this section (Theorem \ref{criterio}) is a consequence of
the following powerful characterization of mixed summable sequences due to
Maurey \cite{Maurey} (see also \cite[16.4.3]{pp1}):

\begin{theorem}
[Maurey]\label{caracter} Let $0<q<s<\infty.$ A sequence $\left(  z_{j}\right)
_{j=1}^{\infty}$ in $E$ is mixed $(s,q)$-summable if, and only if,
\[
\left(  \left(  \int_{B_{E^{\ast}}}\left\vert \left\langle \varphi
,z_{j}\right\rangle \right\vert ^{s}d\mu(\varphi)\right)  ^{\frac{1}{s}%
}\right)  _{j=1}^{\infty}\in\ell_{q}\text{ whenever }\mu\in W(B_{E^{\ast}}).
\]
Besides
\[
\left\Vert \left(  z_{j}\right)  _{j=1}^{\infty}\right\Vert _{mx(s,q)}%
=\sup_{\mu\in W(B_{E^{\ast}})}\left(  \sum_{j=1}^{\infty}\left(
\int_{B_{E^{\ast}}}\left\vert \left\langle \varphi,z_{j}\right\rangle
\right\vert ^{s}d\mu(\varphi)\right)  ^{\frac{q}{s}}\right)  ^{\frac{1}{q}}.
\]
The next theorem shows that our concept has a characterization similar to the
linear case (see \cite{Flore}):
\end{theorem}

\begin{theorem}
\label{criterio}Let $0<p_{1},\ldots,p_{n}\leq q\leq s<\infty.$ An operator
$A\in\mathcal{L}(E_{1},\ldots,E_{n};F)$ is multiple $(s,q;p_{1},\ldots,p_{n})$
mixing summing if, and only if, there is a constant $\sigma\geq0$ such that%
\begin{align}
&  \left(  \sum_{j_{1},\ldots,j_{n}=1}^{m}\left(  \sum_{j=1}^{k}\left\vert
\left\langle \varphi_{j},A(x_{j_{1}}^{(1)},\ldots,x_{j_{n}}^{(n)}%
)\right\rangle \right\vert ^{s}\right)  ^{\frac{q}{s}}\right)  ^{\frac{1}{q}%
}\label{jj}\\
&  \leq\sigma\prod_{l=1}^{n}\left\Vert (x_{i}^{(l)})_{i=1}^{m}\right\Vert
_{w,p_{l}}\left\Vert (\varphi_{j})_{j=1}^{k}\right\Vert _{s}\nonumber
\end{align}
for all $k,m\in\mathbb{N}$, $x_{i}^{(l)}\in E_{l};$ $i=1,\ldots,m,$
$l=1,\ldots,n$ and $\varphi_{j}\in F^{\ast}$ with $j=1,\ldots,k.$ Furthermore,%
\[
\left\Vert A\right\Vert _{mx(s,q;p_{1},\ldots,p_{n})}=\inf\sigma.
\]

\end{theorem}

\begin{proof}
We split the proof into two cases.

(i) Case $s=q.$

From (\ref{jj}) we conclude that
\[
\left(  \sum_{j_{1},\ldots,j_{n}=1}^{m}\left\vert \left\langle \varphi
,A(x_{j_{1}}^{(1)},\ldots,x_{j_{n}}^{(n)})\right\rangle \right\vert
^{q}\right)  ^{\frac{1}{q}}\leq\sigma\prod_{l=1}^{n}\left\Vert (x_{i}%
^{(l)})_{i=1}^{m}\right\Vert _{w,p_{l}}%
\]
for all $\varphi\in B_{F^{\ast}}.$ Thus
\begin{equation}
\left\Vert \left(  A(x_{j_{1}}^{(1)},\ldots,x_{j_{n}}^{(n)})\right)
_{j_{1},...,j_{n}\in\mathbb{N}_{m}}\right\Vert _{w,q}\leq\sigma\prod_{l=1}%
^{n}\left\Vert (x_{i}^{(l)})_{i=1}^{m}\right\Vert _{w,p_{l}} \label{I}%
\end{equation}
and so by Theorem \ref{caracter} and by (\ref{I}) we obtain
\begin{align*}
&  \left\Vert \left(  A(x_{j_{1}}^{(1)},\ldots,x_{j_{n}}^{(n)})\right)
_{j_{1},...,j_{n}\in\mathbb{N}_{m}}\right\Vert _{mx(q,q)}\\
&  =\sup_{\mu\in W(B_{F^{\ast}})}\left(  \sum_{j_{1},\ldots,j_{n}=1}%
^{m}\left(  \int_{B_{F^{\ast}}}\left\vert \left\langle \varphi,A(x_{j_{1}%
}^{(1)},\ldots,x_{j_{n}}^{(n)})\right\rangle \right\vert ^{q}d\mu
(\varphi)\right)  ^{\frac{q}{q}}\right)  ^{\frac{1}{q}}\\
&  \leq\sup_{\mu\in W(B_{F^{\ast}})}\left(  \int_{B_{F^{\ast}}}\sup_{\psi\in
B_{F^{\ast}}}\sum_{j_{1},\ldots,j_{n}=1}^{m}\left\vert \left\langle
\psi,A(x_{j_{1}}^{(1)},\ldots,x_{j_{n}}^{(n)})\right\rangle \right\vert
^{q}d\mu(\varphi)\right)  ^{\frac{1}{q}}\\
&  =\sup_{\mu\in W(B_{F^{\ast}})}\left(  \int_{B_{F^{\ast}}}\left\Vert \left(
A(x_{j_{1}}^{(1)},\ldots,x_{j_{n}}^{(n)})\right)  _{j_{1},...,j_{n}%
\in\mathbb{N}_{m}}\right\Vert _{w,q}^{q}d\mu(\varphi)\right)  ^{\frac{1}{q}}\\
&  \leq\left\Vert \left(  A(x_{j_{1}}^{(1)},\ldots,x_{j_{n}}^{(n)})\right)
_{j_{1},...,j_{n}\in\mathbb{N}_{m}}\right\Vert _{w,q}\\
&  \leq\sigma\prod_{l=1}^{n}\left\Vert (x_{i}^{(l)})_{i=1}^{m}\right\Vert
_{w,p_{l}}.
\end{align*}
Hence, $A\in\Pi_{mx(q,q;p_{1},\ldots,p_{n})}(E_{1},\ldots,E_{n};F)$ and
$\left\Vert A\right\Vert _{mx(q,q;p_{1},\ldots,p_{n})}\leq\sigma$.

Conversely, suppose that $A\in\Pi_{mx(q,q;p_{1},\ldots,p_{n})}(E_{1}%
,\ldots,E_{n};F)$. Given
\[
x_{1}^{(1)},\ldots,x_{m}^{(1)}\in E_{1},\ldots,x_{1}^{(n)},\ldots,x_{m}%
^{(n)}\in E_{n}%
\]
and $\varphi_{1},\ldots,\varphi_{k}\in F^{\ast}$, if
\[
A(x_{j_{1}}^{(1)},\ldots,x_{j_{n}}^{(n)})=\tau_{j_{1},\ldots,j_{n}}%
.y_{j_{1},\ldots,j_{n}},
\]
where $(\tau_{j_{1},\ldots,j_{n}})_{j_{1},...,j_{n}\in\mathbb{N}}\in
\ell_{\infty}$ and $(y_{j_{1},\ldots,j_{n}})_{j_{1},...,j_{n}\in\mathbb{N}}%
\in\ell_{q}^{w}(F;\mathbb{N}^{n})$ we have%
\begin{align*}
&  \left(  \sum_{j_{1},\ldots,j_{n}=1}^{m}\left(  \sum_{j=1}^{k}\left\vert
\left\langle \varphi_{j},A(x_{j_{1}}^{(1)},\ldots,x_{j_{n}}^{(n)}%
)\right\rangle \right\vert ^{q}\right)  ^{\frac{q}{q}}\right)  ^{\frac{1}{q}%
}\\
&  =\left(  \sum_{j=1}^{k}\left(  \left\Vert \varphi_{j}\right\Vert ^{q}%
\sum_{j_{1},\ldots,j_{n}=1}^{m}\left\vert \left\langle \frac{\varphi_{j}%
}{\left\Vert \varphi_{j}\right\Vert },\tau_{j_{1},\ldots,j_{n}}y_{j_{1}%
,\ldots,j_{n}}\right\rangle \right\vert ^{q}\right)  \right)  ^{\frac{1}{q}}\\
&  =\left(  \sum_{j=1}^{k}\left\Vert \varphi_{j}\right\Vert ^{q}\right)
^{\frac{1}{q}}\left(  \sum_{j_{1},\ldots,j_{n}=1}^{m}\left\vert \tau
_{j_{1},\ldots,j_{n}}\right\vert ^{q}\left\vert \left\langle \frac{\varphi
_{j}}{\left\Vert \varphi_{j}\right\Vert },y_{j_{1},\ldots,j_{n}}\right\rangle
\right\vert ^{q}\right)  ^{\frac{1}{q}}\\
&  \leq\left\Vert (\varphi_{j})_{j=1}^{k}\right\Vert _{q}\left\Vert
(\tau_{j_{1},\ldots,j_{n}})_{j\in\mathbb{N}^{n}}\right\Vert _{\infty
}\left\Vert (y_{j_{1},\ldots,j_{n}})_{j\in\mathbb{N}^{n}}\right\Vert _{w,q}.
\end{align*}
Taking the infimum in both sides, we obtain
\begin{align*}
&  \left(  \sum_{j_{1},\ldots,j_{n}=1}^{m}\left(  \sum_{j=1}^{k}\left\vert
\left\langle \varphi_{j},A(x_{j_{1}}^{(1)},\ldots,x_{j_{n}}^{(n)}%
)\right\rangle \right\vert ^{q}\right)  ^{\frac{q}{q}}\right)  ^{\frac{1}{q}%
}\\
&  \leq\left\Vert (\varphi_{j})_{j=1}^{k}\right\Vert _{q}\left\Vert \left(
A(x_{j_{1}}^{(1)},\ldots,x_{j_{n}}^{(n)})\right)  _{j\in\mathbb{N}_{m}^{n}%
}\right\Vert _{m,(q,q)}\\
&  \leq\left\Vert (\varphi_{j})_{j=1}^{k}\right\Vert _{q}\left\Vert
A\right\Vert _{mx(q,q;p_{1},\ldots,p_{n})}\prod_{l=1}^{n}\left\Vert
(x_{i}^{(l)})_{i=1}^{m}\right\Vert _{w,p_{l}}.
\end{align*}
Therefore $\inf\sigma\leq\left\Vert A\right\Vert _{mx(q,q;p_{1},\ldots,p_{n}%
)}$ and with the last inequality we obtain
\[
\left\Vert A\right\Vert _{mx(q,q;p_{1},\ldots,p_{n})}=\inf\sigma.
\]

(ii) Case $s>q.$

Let $A\in\Pi_{mx(s,q;p_{1},\ldots,p_{n})}(E_{1},\ldots,E_{n};F).$ Given
$0\neq\varphi_{1},\ldots,\varphi_{k}\in F^{\ast}$ we define the probability
measure%
\[
\nu=\sum_{j=1}^{k}\nu_{j}\delta_{j},\text{ where }\nu_{j}=\frac{\left\Vert
\varphi_{j}\right\Vert ^{s}}{\sum_{j=1}^{k}\left\Vert \varphi_{j}\right\Vert
^{s}}%
\]
and $\delta_{j}$ is the Dirac measure at the point $\tilde{\varphi}_{j}%
=\frac{\varphi_{j}}{\left\Vert \varphi_{j}\right\Vert }.$

For $x_{1}^{(1)},\ldots,x_{m}^{(1)}\in E_{1},\ldots$, $x_{1}^{(n)}%
,\ldots,x_{m}^{(n)}\in E_{n}$, note that
\begin{align*}
&  \int_{B_{F^{\ast}}}\left\vert \left\langle \varphi,A(x_{j_{1}}^{(1)}%
,\ldots,x_{j_{n}}^{(n)})\right\rangle \right\vert ^{s}d\nu(\varphi)\\
&  =\sum_{j=1}^{k}\left\vert \left\langle \tilde{\varphi_{j}},A(x_{j_{1}%
}^{(1)},\ldots,x_{j_{n}}^{(n)})\right\rangle \right\vert ^{s}\nu
(\tilde{\varphi_{j}})\\
&  =\sum_{j=1}^{k}\left\vert \left\langle \frac{\varphi_{j}}{\left\Vert
\varphi_{j}\right\Vert },A(x_{j_{1}}^{(1)},\ldots,x_{j_{n}}^{(n)}%
)\right\rangle \right\vert ^{s}.\nu_{j}.\delta_{j}(\tilde{\varphi_{j}})\\
&  =\sum_{j=1}^{k}\left\vert \left\langle \frac{\varphi_{j}}{\left\Vert
\varphi_{j}\right\Vert },A(x_{j_{1}}^{(1)},\ldots,x_{j_{n}}^{(n)}%
)\right\rangle \right\vert ^{s}.\frac{\left\Vert \varphi_{j}\right\Vert ^{s}%
}{\sum_{j=1}^{k}\left\Vert \varphi_{j}\right\Vert ^{s}}\\
&  =\frac{1}{\left\Vert (\varphi_{j})_{j=1}^{k}\right\Vert _{s}^{s}}\sum
_{j=1}^{k}\left\vert \left\langle \varphi_{j},A(x_{j_{1}}^{(1)},\ldots
,x_{j_{n}}^{(n)})\right\rangle \right\vert ^{s}.
\end{align*}

From the previous equalities and from Theorem \ref{caracter} we have
\begin{align*}
&  \left(  \sum_{j_{1},\ldots,j_{n}=1}^{m}\left(  \sum_{j=1}^{k}\left\vert
\left\langle \varphi_{j},A(x_{j_{1}}^{(1)},\ldots,x_{j_{n}}^{(n)}%
)\right\rangle \right\vert ^{s}\right)  ^{\frac{q}{s}}\right)  ^{\frac{1}{q}%
}\\
&  =\left(  \sum_{j_{1},\ldots,j_{n}=1}^{m}\left(  \int_{B_{F^{\ast}}%
}\left\vert \left\langle \varphi,A(x_{j_{1}}^{(1)},\ldots,x_{j_{n}}%
^{(n)})\right\rangle \right\vert ^{s}d\nu(\varphi)\right)  ^{\frac{q}{s}%
}\right)  ^{\frac{1}{q}}\left\Vert (\varphi_{j})_{j=1}^{k}\right\Vert _{s}\\
&  \leq\left\Vert \left(  A(x_{j_{1}}^{(1)},\ldots,x_{j_{n}}^{(n)})\right)
_{j\in\mathbb{N}_{m}^{n}}\right\Vert _{m,(s,q)}\left\Vert (\varphi_{j}%
)_{j=1}^{k}\right\Vert _{s}\\
&  \leq\left\Vert A\right\Vert _{mx(s,q;p_{1},\ldots,p_{n})}\prod_{l=1}%
^{n}\left\Vert (x_{i}^{(l)})_{i=1}^{m}\right\Vert _{w,p_{l}}\left\Vert
(\varphi_{j})_{j=1}^{k}\right\Vert _{s}.
\end{align*}
and we obtain (\ref{jj}) with $\inf\sigma\leq\left\Vert A\right\Vert
_{mx(s,q;p_{1},\ldots,p_{n})}.$

Reciprocally, with the same idea and using (\ref{jj}), given $\nu=\sum
_{i=1}^{k}\nu_{i}\delta_{i}$ a discrete probability measure onto $B_{F^{\ast}%
}$ we obtain
\begin{align*}
&  \left(  \sum_{j_{1},\ldots,j_{n}=1}^{m}\left(  \int_{B_{F^{\ast}}%
}\left\vert \left\langle \varphi,A(x_{j_{1}}^{(1)},\ldots,x_{j_{n}}%
^{(n)})\right\rangle \right\vert ^{s}d\nu(\varphi)\right)  ^{\frac{q}{s}%
}\right)  ^{\frac{1}{q}}\\
&  =\left(  \sum_{j_{1},\ldots,j_{n}=1}^{m}\left(  \sum_{j=1}^{k}\left\vert
\left\langle \nu_{j}^{\frac{1}{s}}\varphi_{j},A(x_{j_{1}}^{(1)},\ldots
,x_{j_{n}}^{(n)})\right\rangle \right\vert ^{s}\right)  ^{\frac{q}{s}}\right)
^{\frac{1}{q}}\\
&  \leq\sigma\prod_{l=1}^{n}\left\Vert (x_{i}^{(l)})_{i=1}^{m}\right\Vert
_{w,p_{l}}\left\Vert (\nu_{j}^{\frac{1}{s}}\varphi_{j})_{j=1}^{k}\right\Vert
_{s}\\
&  \leq\sigma\prod_{l=1}^{n}\left\Vert (x_{i}^{(l)})_{i=1}^{m}\right\Vert
_{w,p_{l}}.
\end{align*}
The previous inequality holds for every $\nu\in W(B_{F^{\ast}})$, since the
discrete probability measures are dense in $W(B_{F^{\ast}})$. Therefore, from
Theorem \ref{caracter} we obtain
\[
\left\Vert \left(  A(x_{j_{1}}^{(1)},\ldots,x_{j_{n}}^{(n)})\right)
_{j\in\mathbb{N}_{m}^{n}}\right\Vert _{mx(s,q)}\leq\sigma\prod_{l=1}%
^{n}\left\Vert (x_{i}^{(l)})_{i=1}^{m}\right\Vert _{w,p_{l}},
\]
for all $m\in\mathbb{N}$ and%
\[
\left\Vert A\right\Vert _{mx(s,q;p_{1},\ldots,p_{n})}=\inf\sigma.
\]

\end{proof}

\subsection{A quotient theorem}

For linear operators, $S\in\mathcal{L}\left(  E;F\right)  $ is $\left(
s,p\right)  $-mixing summing if and only if $TS$ is absolutely $p$-summing for
all $T\in\Pi_{s}\left(  F;G\right)  $. In other words%
\[
\Pi_{mx\left(  s,p\right)  }\left(  E;F\right)  =\left(  \Pi_{s}\left(
F;G\right)  \right)  ^{-1}\circ\Pi_{p}\left(  E;G\right)  .
\]

For details we refer to \cite[Section 32]{Flore} and \cite{pp1}. In this
section we show that our approach provides a perfect multilinear extension of
this result. We show that the following assertions are equivalent:

\begin{itemize}
\item $T\in\mathcal{L}(E_{1},\ldots,E_{n};F)$ is multiple $\left(  s,q;p_{1}
,\ldots,p_{n}\right)  $-mixing summing.

\item $u\circ T\in\mathcal{L}_{mas(q;p_{1},\ldots,p_{n})}(E_{1},\ldots
,E_{n};G)$ for all $u\in\Pi_{s}\left(  F;G\right)  $ and $T\in\mathcal{L}%
(E_{1},\ldots,E_{n};F).$
\end{itemize}

Using a different notation, we will show the following quotient theorem:
\begin{equation}
\Pi_{mx\left(  s,q;p_{1},\ldots,p_{n}\right)  }\left(  E_{1},\ldots
,E_{n};F\right)  =\left(  \Pi_{s}\left(  F;G\right)  \right)  ^{-1}%
\circ\mathcal{L}_{mas\left(  q;p_{1},\ldots,p_{n}\right)  }\left(
E_{1},\ldots,E_{n};G\right)  \label{dtt}%
\end{equation}
for all $E_{1},\ldots,E_{n},F$ and $G.$

The quotient theorem (\ref{dtt}) is a direct consequence of the forthcoming
Propositions \ref{pp9} and \ref{pp99}. First we need the following lemma:

\begin{lemma}
\label{jlk}Let $A\in\mathcal{L}(E_{1},\ldots,E_{n};F)$ be so that
\[
u\circ A\in\mathcal{L}_{mas(p;p_{1},\ldots,p_{n})}(E_{1},\ldots,E_{n};G)
\]
for all Banach space $G$ and all $u\in\Pi_{r}(F;G).$ Then, there is a $C\geq0$
such that%
\begin{equation}
\left\Vert u\circ A\right\Vert _{(p;p_{1},\ldots,p_{n})}\leq C\pi
_{r}(u).\label{gtr}%
\end{equation}

\end{lemma}

\begin{proof}
Suppose that (\ref{gtr}) is not true. So, for all positive integer $k$ there
exist Banach spaces $F_{k}$ and $u_{k}\in\Pi_{r}(F;F_{k})$ so that%
\[
\pi_{r}(u_{k})\leq\frac{1}{2^{k}}\text{ and }\left\Vert u_{k}\circ
A\right\Vert _{(p;p_{1},\ldots,p_{n})}\geq k.
\]
Let $J_{k}:F_{k}\rightarrow\ell_{2}\left(  \left(  F_{k}\right)
_{k=1}^{\infty}\right)  $ and $Q_{j}:\ell_{2}\left(  \left(  F_{k}\right)
_{k=1}^{\infty}\right)  \rightarrow F_{j}$ be the canonical maps for all
positive integers $j,k$. Since%
\[
\pi_{r}\left(  \sum\limits_{k=n_{1}}^{n_{2}}J_{k}\circ u_{k}\right)  \leq
\sum\limits_{k=n_{1}}^{n_{2}}\pi_{r}\left(  J_{k}\circ u_{k}\right)  \leq
\sum\limits_{k=n_{1}}^{n_{2}}\pi_{r}\left(  u_{k}\right)  \leq\sum
\limits_{k=n_{1}}^{n_{2}}\frac{1}{2^{k}}%
\]
it follows that
\[
u:=\sum\limits_{j=1}^{\infty}J_{j}\circ u_{j}\in\Pi_{r}(F;\ell_{2}\left(
\left(  F_{k}\right)  _{k=1}^{\infty}\right)  ).
\]
Since $u_{k}=Q_{k}\circ u$, we thus have%
\[
k\leq\left\Vert u_{k}\circ A\right\Vert _{(p;p_{1},\ldots,p_{n})}=\left\Vert
Q_{k}\circ u\circ A\right\Vert _{(p;p_{1},\ldots,p_{n})}\leq\left\Vert u\circ
A\right\Vert _{(p;p_{1},\ldots,p_{n})},
\]
a contradiction.
\end{proof}

\begin{proposition}
\label{pp9}If $A\in\mathcal{L}(E_{1},\ldots,E_{n};F)$ is so that $u\circ
A\in\mathcal{L}_{mas(q;p_{1},\ldots.,p_{n})}(E_{1},\ldots,E_{n};G)$ for all
$u\in\Pi_{s}\left(  F;G\right)  ,$ then%
\[
A\in\Pi_{mx(s,q;p_{1},\ldots.,p_{n})}(E_{1},\ldots,E_{n};F).
\]

\end{proposition}

\begin{proof}
Let $x_{i}^{(j)}\in E_{j}$ with $\left(  i,j\right)  \in\left\{
1,\ldots,m\right\}  \times\left\{  1,\ldots,n\right\}  .$ Consider
$S:F\rightarrow\ell_{s}^{k}$ defined by
\[
S(y)=\left(  \varphi_{j}\left(  y\right)  \right)  _{j=1}^{k}.
\]
It is not difficult to show that
\[
\pi_{s}\left(  S\right)  \leq\left\Vert \left(  \varphi_{j}\right)  _{j=1}%
^{k}\right\Vert _{s}.
\]
Since $S\circ A\in\mathcal{L}_{mas(q;p_{1},\ldots.,p_{n})}(E_{1},\ldots
,E_{n};\ell_{s}^{k})$ and invoking Lemma \ref{jlk}, there is a constant $C>0$
so that%
\begin{align*}
&  \left(  \sum_{j_{1},\ldots,j_{n}=1}^{m}\left(  \sum_{j=1}^{k}\left\vert
\left\langle \varphi_{j},A(x_{j_{1}}^{(1)},\ldots,x_{j_{n}}^{(n)}%
)\right\rangle \right\vert ^{s}\right)  ^{\frac{q}{s}}\right)  ^{\frac{1}{q}%
}\\
&  =\left(  \sum_{j_{1},\ldots,j_{n}=1}^{m}\left\Vert S\circ A\left(
x_{j_{1}}^{(1)},\ldots,x_{j_{n}}^{(n)}\right)  \right\Vert _{s}^{q}\right)
^{\frac{1}{q}}\\
&  \leq\left\Vert S\circ A\right\Vert _{(q;p_{1},\ldots.,p_{n})}%
\prod\limits_{j=1}^{n}\left\Vert \left(  x_{i}^{(j)}\right)  _{i=1}%
^{m}\right\Vert _{w,p_{j}}\\
&  \leq C\pi_{s}\left(  S\right)  \prod\limits_{j=1}^{n}\left\Vert \left(
x_{i}^{(j)}\right)  _{i=1}^{m}\right\Vert _{w,p_{j}}\\
&  \leq C\left\Vert \left(  \varphi_{j}\right)  _{j=1}^{k}\right\Vert
_{s}\prod\limits_{j=1}^{n}\left\Vert \left(  x_{i}^{(j)}\right)  _{i=1}%
^{m}\right\Vert _{w,p_{j}}.
\end{align*}

\end{proof}

\begin{proposition}
\label{pp99}If $A\in\Pi_{mx(s,q;p_{1},\ldots.,p_{n})}(E_{1},\ldots,E_{n};F),$
then%
\begin{equation}
u\circ A\in\Pi_{(q;p_{1},\ldots.,p_{n})}(E_{1},\ldots,E_{n};G) \label{ggr}%
\end{equation}
and
\begin{equation}
\left\Vert u\circ A\right\Vert _{(q;p_{1},\ldots.,p_{n})}\leq\pi_{s}\left(
u\right)  \left\Vert A\right\Vert _{mx(s,q;p_{1},\ldots.,p_{n})} \label{ggs}%
\end{equation}
for all $u\in\Pi_{s}\left(  F;G\right)  .$
\end{proposition}

\begin{proof}
Let $x_{i}^{(j)}\in E_{j}$ with $\left(  i,j\right)  \in\left\{
1,\ldots,m\right\}  \times\left\{  1,\ldots,n\right\}  .$ Given $\varepsilon
>0$ there are $\tau_{j_{1},\ldots,j_{n}}\in K$ and $y_{j_{1},\ldots,j_{n}}\in
F$ so that%
\[
A(x_{j_{1}}^{(1)},\ldots,x_{j_{n}}^{(n)})=\tau_{j_{1},\ldots,j_{n}}%
y_{j_{1},\ldots,j_{n}}%
\]
and%
\begin{align*}
&  \left\Vert \left(  \tau_{j_{1},\ldots,j_{n}}\right)  _{j_{1},\ldots
,j_{n}=1}^{m}\right\Vert _{r}\left\Vert \left(  y_{j_{1},\ldots,j_{n}}\right)
_{j_{1},\ldots,j_{n}=1}^{m}\right\Vert _{w,s}\\
&  <\left(  1+\varepsilon\right)  \left\Vert \left(  A(x_{j_{1}}^{(1)}%
,\ldots,x_{j_{n}}^{(n)})\right)  _{j_{1},\ldots,j_{n}=1}^{m}\right\Vert
_{mx\left(  s,q)\right)  }\\
&  \leq\left(  1+\varepsilon\right)  \left\Vert A\right\Vert _{mx(s,q;p_{1}%
,\ldots.,p_{n})}\prod\limits_{j=1}^{n}\left\Vert \left(  x_{i}^{(j)}\right)
_{i=1}^{m}\right\Vert _{w,p_{j}}.
\end{align*}
Hence, using H\"{o}lder's Inequality we obtain
\begin{align*}
&  \left\Vert \left(  u\circ A(x_{j_{1}}^{(1)},\ldots,x_{j_{n}}^{(n)})\right)
_{j_{1},\ldots,j_{n}=1}^{m}\right\Vert _{q}\\
&  \leq\left\Vert \left(  \tau_{j_{1},\ldots,j_{n}}\right)  _{j_{1}%
,\ldots,j_{n}=1}^{m}\right\Vert _{r}\left\Vert \left(  u\left(  y_{j_{1}%
,\ldots,j_{n}}\right)  \right)  _{j_{1},\ldots,j_{n}=1}^{m}\right\Vert _{s}\\
&  \leq\left(  1+\varepsilon\right)  \pi_{s}\left(  u\right)  \left\Vert
A\right\Vert _{mx(s,q;p_{1},\ldots.,p_{n})}\prod\limits_{j=1}^{n}\left\Vert
\left(  x_{i}^{(j)}\right)  _{i=1}^{m}\right\Vert _{w,p_{j}}%
\end{align*}
and making $\varepsilon\rightarrow0$ we get (\ref{ggr}) and (\ref{ggs}).
\end{proof}

\subsection{Coherence and compatibility}

The polynomial version of multiple mixing summing operators can be stated by
using the symmetric multilinear operator associated to the polynomials:

\begin{definition}
Let $0<p\leq s<\infty.$ A polynomial $P\in\mathcal{P}(^{n}E;F)$ is multiple
$(s,p)$-mixing summing if $\check{P}$ is multiple $(s,p;p)$-mixing summing.
Besides,%
\[
\left\Vert P\right\Vert _{mx(s,p)}:=\left\Vert \check{P}\right\Vert
_{mx(s,p;p)}.
\]

\end{definition}

The following proposition, whose proof is standard, shows that, as it happens
to multiple summing multilinear operators, coincidence results for multiple
mixing summing multilinear operators imply in coincidence results for smaller degrees:

\begin{proposition}
If $\mathcal{L}(E_{1},\ldots,E_{n};F)=\Pi_{mx(s,q;p_{1},\ldots,p_{n})}%
(E_{1},\ldots,E_{n};F)$ then%
\[
\mathcal{L}(E_{k_{1}},\ldots,E_{k_{j}};F)=\Pi_{mx(s,q;p_{k_{1}},\ldots
,p_{k_{j}})}(E_{k_{1}},\ldots,E_{k_{j}};F)
\]
whenever $1\leq j<n$ and $\left\{  k_{1}<\cdots<k_{j}\right\}  \subset\left\{
1,\ldots,n\right\}  $.
\end{proposition}

Similarly to the previous section one can show that $\left(  \mathcal{P}%
_{mx(s,p)}^{n},\left\Vert \cdot\right\Vert _{mx(s,p)}\right)  _{n=1}^{\infty}$
is coherent and for each $n$ it is compatible with the operator ideal $\left(
\Pi_{mx(s,p)},\pi_{mx(s,p)}\right)  $. For example we prove (i) of Definition
\ref{IdeaisCoerentes}:

\begin{proposition}
If $P\in\mathcal{P}_{mx(s,p)}(^{n}E;F)$ and $a\in E$, then $P_{a}%
\in\mathcal{P}_{mx(s,p)}(^{n-1}E;F)$ and
\[
\left\Vert P_{a}\right\Vert _{mx(s,p)}\leq\left\Vert P\right\Vert
_{mx(s,p)}\left\Vert a\right\Vert .
\]

\end{proposition}

\begin{proof}
Since $\check{P}\in\Pi_{mx(s,p)}(^{n}E;F)$ we have%
\[
\left(  \sum_{j_{1},\ldots,j_{n}=1}^{m}\left(  \sum_{j=1}^{k}\left\vert
\left\langle \varphi_{j},\check{P}(x_{j_{1}}^{(1)},\ldots,x_{j_{n}}%
^{(n)})\right\rangle \right\vert ^{s}\right)  ^{\frac{p}{s}}\right)
^{\frac{1}{p}}\leq\sigma\prod_{l=1}^{n}\left\Vert (x_{i}^{(l)})_{i=1}%
^{m}\right\Vert _{w,p}\left\Vert (\varphi_{j})_{j=1}^{k}\right\Vert _{s}.
\]
and by choosing $x_{1}^{(n)}=a$ and $x_{j}^{(n)}=0$ for $j>1$ we have%
\begin{align*}
&  \left(  \sum_{j_{1},\ldots,j_{n-1}=1}^{m}\left(  \sum_{j=1}^{k}\left\vert
\left\langle \varphi_{j},\check{P}_{a}(x_{j_{1}}^{(1)},\ldots,x_{j_{n-1}%
}^{(n-1)})\right\rangle \right\vert ^{s}\right)  ^{\frac{p}{s}}\right)
^{\frac{1}{p}}\\
&  =\left(  \sum_{j_{1},\ldots,j_{n-1}=1}^{m}\left(  \sum_{j=1}^{k}\left\vert
\left\langle \varphi_{j},\check{P}(x_{j_{1}}^{(1)},\ldots,x_{j_{n-1}}%
^{(n-1)},a)\right\rangle \right\vert ^{s}\right)  ^{\frac{p}{s}}\right)
^{\frac{1}{p}}\\
&  =\left(  \sum_{j_{1},\ldots,j_{n}=1}^{m}\left(  \sum_{j=1}^{k}\left\vert
\left\langle \varphi_{j},\check{P}(x_{j_{1}}^{(1)},\ldots,x_{j_{n}}%
^{(n)})\right\rangle \right\vert ^{s}\right)  ^{\frac{p}{s}}\right)
^{\frac{1}{p}}\\
&  \leq\left\Vert P\right\Vert _{mx(s,p)}\left\Vert a\right\Vert \prod
_{l=1}^{n-1}\left\Vert (x_{i}^{(l)})_{i=1}^{m}\right\Vert _{w,p}\left\Vert
(\varphi_{j})_{j=1}^{k}\right\Vert _{s}.
\end{align*}

\end{proof}

\section{Final comments and directions for further research}

The concepts of multiple mixing summing and multiple $\left(  p;q;r_{1}%
,\ldots,r_{n}\right)  $-summing polynomials/multilinear operators, as natural
extensions of the notion of multiple summing multilinear operators, can be
further investigated following different directions: coincidence theorems,
generalizations to holomorphic mappings, or inclusion theorems, among others.

The study of coincidence theorems may follow the lines of \cite{Na} combined
with the results from the respective linear theories; the study of holomorphic
mappings may follow \cite{Junek} and for inclusion theorems \cite{davidstudia}
is certainly a good source of inspiration.

We encourage the interested reader to investigate other variants of mixing
summability and $\left(  p;q;r_{1},\ldots,r_{n}\right)  $-summability
following the lines given in \cite{jgdd, Nach00, MST, parc}.

\end{document}